\documentclass[12pt]{imsart}
\RequirePackage{natbib}
\usepackage{bbm}
\usepackage{amsmath,amssymb,amsthm}
\usepackage{graphics,epsfig}
\usepackage{hyperref}
\usepackage{natbib}
\usepackage{color}
\usepackage{graphicx}
\usepackage{caption}
\usepackage{subcaption}
\usepackage{float}
\usepackage{dsfont}
\usepackage{mathrsfs}
\usepackage{multirow}
\usepackage{comment}
\usepackage{bm}
\usepackage{geometry}

\def \qq {\boldsymbol{q}}

\def \ra {\rightarrow}

\def \E {\mathbb{E}}

\def \a {\alpha}
\def \be {\beta}

\newtheorem{definition}{\bf Definition}

	\newtheorem{remark}{\bf Remark}

	\newtheorem{theorem}{\bf Theorem}
	
	\newtheorem{prop}{\bf Proposition}
	
	\newtheorem{as}{\bf Assumption}

\renewcommand{\epsilon}{\varepsilon}

\begin{document}

\begin{frontmatter}

\title{Stochastic Control of a SIR Model with Non-linear Incidence Rate Through Euclidean Path Integral}
\runtitle{Pandemic control}

\begin{aug}
\author{\fnms{Paramahansa} 
	\snm{Pramanik}
	\ead[label=e1]{ppramanik@southalabama.edu}}

\runauthor{P. Pramanik}

\affiliation{University of South Alabama}

\address{Department of Mathematics and Statistics\\ University of South Alabama\\ Mobile, AL 36688 USA.\\email:ppramanik@southalabama.edu\\Phone:251-341-3098.}

\end{aug}

\begin{abstract}
	This paper utilizes a stochastic Susceptible-Infected-recovered (SIR) model with a non-linear incidence rate to perform a detailed mathematical study of optimal lock-down intensity and vaccination rate under the COVID-19 pandemic environment. We use a Feynman-type path integral control approach to determine a Fokker-Plank type equation of this system. Since we assume the availability of information on the COVID-19 pandemic is complete and perfect, we can show a unique fixed point. A non-linear incidence rate is used because, it can be raised from saturation effects that if the proportion of infected agents is very high so that exposure to the pandemic is inevitable, then the transmission rate responds slower than linearity to the increase in the number of infections. The simulation study shows that with higher diffusion coefficients susceptible and recovery curves keep the downward trends while the infection curve becomes ergodic. Finally, we perform a data analysis using UK data at the beginning of 2021 and compare it with our theoretical results. 
\end{abstract}



\end{frontmatter}

\section{Introduction}
In current days we see \emph{locking downs} of economies and increasing \emph{vaccination rate} as a strategy to reduce the spread of COVID-19 which already has claimed more than one million lives in the United States and more than six million across the globe. During the past couple of years, it has been clear that an adequate synthesis requires better epidemiology, better economic analysis, and more advanced optimization techniques to tame this pandemic. Almost all mathematical methods of epidemic models descend from the susceptible-infected-recovered (SIR) model \citep{kermack1927, rao2014}. The dynamic behavior of different epidemic models has been studied extensively \citep{acemoglu2020,caulkins2021optimal,pramanik2022lock,rao2014}. Stochastic pandemic modeling is important when the number of infected agents is small or when the different transmission and recovery rates influence the pandemic outcome \citep{allen2017}.

In this paper, we perform a Feynman-type path integral approach for a recursive formulation of a quadratic cost function with a \emph{forward-looking} stochastic SIR model and an infection dynamics based on an Erdos-Renyi type random network. A  Fokker-Plank type  equation is obtained for this COVID-19 environment which is analogous to an HJB equation \citep{yeung2006} and a saddle-point functional equation \citep{marcet2019}. Solving for the first order condition of the Fokker-Plank equation the optimal lock-down intensity and vaccination rate are obtained.  As the movement of infection in a society over time is stochastic, analogous to the movement of a quantum particle, a Feynman-type path integral method of quantum physics has been used to determine optimal \emph{lock-down intensity} and \emph{vaccination rate}. We define \emph{lock-down intensity} as the ratio of employment due to COVID-19 to total employment under the absence of the pandemic. Therefore, the value of this lock-down intensity lies between $0$ and $1$ where $0$ stands for complete shut-down of an economy. Our formulation is based on path integral control and dynamic programming tools facilitate the analysis and permit the application of an algorithm to obtain a numerical solution for this stochastic pandemic control model. Throughout this paper, we assume all agents in the pandemic environment are risk averse. Therefore, the simulation of optimal lock-down intensity goes up at the beginning of our time interval and then comes very close to zero. The reason behind this is that, due to the availability of \emph{perfect} and \emph{complete} information an individual does not want to go out and get infected by COVID-19.

We have used a non-linear incidence rate because, it can be raised from saturation effects that if the proportion of infected agents is very high so that exposure to the pandemic is inevitable, then the transmission rate responds slower than linear to the increase in the number of infections \citep{korobeinikov2005,rao2014}. In \cite{capasso1978} a \emph{saturated transmission rate} is defined as $b(S,I)=\be SI/(1+\rho I)$, for all proportionality constant $\rho\in(0,1]$, stochastic infection rate $\be\in[0,1]$, where $\be I$ is a measure of pandemic infection force and $1/(1+\rho I)$ is a measure of inhibition effect from the behavioral change of the susceptible agents when their number increases \citep{rao2014}. To become feasible in the biological sense, for all $S, I>0$ assume the function $b(s,I)$ is smooth and concave with respect to $I$ such that, $b(S,0)=b(0,I)=0$, $\partial b/\partial S=\be I/(1+\rho I)>0$, $\partial b/\partial I=\be S/(1+\rho I)^2>0$ and $\partial^2b/\partial I^2=-2\rho\be S/(1+\rho I)^3<0$. The second order condition implies that when the number of infections is very high that the exposure to the pandemic is certain, the incidence rate responds slower than linearity in $I$ \citep{rao2014}. To the best of our knowledge, less amount of research has been done with stochastic perturbation on a SIR  pandemic COVID-19 model with $\be SI/(1+\rho I)$ as a saturated transmission rate.

Feynman path integral is a method of quantization uses a quantum Lagrangian function, while Schr\"odinger quantization uses a Hamiltonian function \citep{fujiwara2017}.  Since the path integral approach provides a different viewpoint from Schr\"odinger's quantization, it is a very useful tool not only in quantum physics but also in engineering, biophysics, economics, and finance \citep{anderson2011,fujiwara2017,kappen2005,yang2014path}. These two approaches are believed to be equivalent but, this equivalence has not fully proved mathematically as the mathematical difficulties lie in the fact that the Feynman path integral is not an integral utilizing a countably additive measure \citep{johnson2000,fujiwara2017,pr2,pr3}. As the complexity and memory requirements of grid-based partial differential equation (PDE) solvers increase exponentially as the dimension of the system increases, this method becomes impractical in the case with high dimensions \citep{pr0,pr1,yang2014path}. As an alternative one can use a Monte Carlo scheme and this is the main idea of path integral control \citep{kappen2005,theodorou2010,theodorou2011,morzfeld2015}. Path integral control solves a class of stochastic control problems with a Monte Carlo method \citep{uddin2020,islam2021,alam2021,alam2022} for an HJB equation and this approach avoids the need for a global grid of the domain of the HJB equation \citep{yang2014path}. If the objective function is quadratic and the differential equations are linear, then the solution is given in terms of several Ricatti equations that can be solved efficiently \citep{kappen2007b,pramanik2020motivation,pramanik2021,pramanik2021scoring}.

Although incorporating randomness with its HJB equation is straightforward, difficulties come due to dimensionality when a numerical solution is calculated for both deterministic and stochastic HJB \citep{kappen2007b}. The general stochastic control problem is intractable to solve computationally as it requires an exponential amount of memory and computational time because, the state space needs to be discretized and hence, becomes exponentially large in the number of dimensions \citep{theodorou2010,theodorou2011,yang2014path}. Therefore, in order to determine the expected values, it is necessary to visit all states which leads to the inefficient summations of exponentially large sums \citep{kappen2007b,yang2014path,pramanik2021}. This is the main reason to implement a path integral control approach to deal with stochastic pandemic control.

Following is the structure of this paper. Section 2.1 describes the main problem formulation. We discuss the properties of stochastic SIR and the quadratic cost function. We also show that under perfect and complete information our SIR model has a unique solution. Section 2.2 discusses about the transmission of COVID-19 in a community with Erdos-Renyi random interaction based on five different immune groups. In this section, we also consider the fine particulate matter to observe the effect of air pollution on an individual infected by the pandemic. Section 2.3 constructs the system of stochastic constraints including infection dynamics and its properties. Section 2.4 describes the main theoretical results of this paper. We did some simulation studies and real data analysis based on UK data for SIR in section 3 based on the results obtained from section 2.4. Finally, section 4 concludes the paper. All the proofs are in the appendix.

\section{Framework}

\medskip

\subsection{Model}

\medskip

In this section, we are going to construct a dynamic framework where a social planner's cost is minimized subject to a stochastic SIR model with pandemic spread dynamics. Throughout the paper, we are considering the stochastic optimization problem of a single agent, and to make our model simple we assume all agents' objectives are identical. Therefore, we ignore any subscripts to represent an agent. Following \cite{lesniewski2020} an agent's objective is to minimize a cost function 

\begin{multline}\label{0}
u^*=\min_{v,e\in\mathcal U}\E_0\biggr\{\int_0^t\left[\exp(-rs)\left[S(s)\left(\mbox{$\frac{1}{2}$}\a_{11}v^2(s)+\a_{12}v(s)+\a_{13}\right)\right.\right.\\
\left.\left.+I(s)\left(\mbox{$\frac{1}{2}$}\a_{21}e^2(s)+\a_{22}e(s)+\a_{23}\right)\right]+\be(e(s),v(s))S(s)I(s)\right]ds\bigg|\mathcal F_0\biggr\},
\end{multline}
subject to
\begin{align}\label{1}
dS(s) & = \biggr\{\eta N(s)-\be(e(s),v(s))\frac{S(s)I(s)}{\left[1+\rho I(s)\right]+\eta N(s)}-\kappa S(s)-v(s)\notag\\&\hspace{1cm}+\zeta R(s)\biggr\}ds+\sigma_1\left[S(s)-S^*\right]dB_1(s),\notag\\ dI(s)&=\left\{\be(e(s),v(s))\frac{S(s)I(s)}{\left[1+\rho I(s)\right]+\eta N(s)}-(\mu+\kappa)I(s)-e(s)\right\}ds\notag\\ &\hspace{1cm}+\sigma_2\left[I(s)-I^*\right]dB_2(s),\notag\\ dR(s) &=\left\{\mu v(s)I(s)-[\kappa +\zeta]e(s) R(s)\right\}ds+\sigma_3\left[R(s)-R^*\right]dB_3(s),
\end{align}
with the stochastic differential infection rate $\be$, a function of vaccination rate $v$ and lock down intensity $e$. In Equation (\ref{0}), $r\in(0,1)$ is a continuous discounting factor; $S$ and $I$ represent percentage of total population ($N$) susceptible to and infected with  COVID-19. $R$ is the percentage of people removed from $N$ where $S+I+R=N$. $R$ includes people who got completely recovered from COVID-19 and those people who passed away because of this pandemic. As $S$, $I$ and $R$ are represented in terms of percentages therefore, $N=100$. Furthermore, in Equation (\ref{0}) $u^*=(v^*,e^*)$ represents an optimal level of vaccination rate and lock-down intensity respectively. The coefficients $\a_{ij}$ for all $i=1,2$ and $j=1,2,3$ are determined by the overall cost functions with $\a_{11}>0$ \citep{rao2014}. Finally, $\mathcal F_0$ is the filtration process of the state variables $S,I$ and $R$ starting at time $0\in[0,t]$. Hence, $\E_0[.]=\E[.|S(0),I(0),R(0);\mathcal F_0]$ where $S(0)$, $I(0)$ and $R(0)$ are the initial conditions.

\medskip

In the System of Equations (\ref{1}) $\eta$ is the birth rate, $1/\left[1+\rho I(s)\right]$ is a measure of inhibition effect from the behavioral change of the susceptible individual, $\kappa$ is the natural death rate, $\zeta$ is the rate at which a recovered person loses immunity and returns to the susceptible class and $\mu$ is the natural recovery rate. $\sigma_1$, $\sigma_2$ and $\sigma_3$ are assumed to be real constants and are defined as the intensity of stochastic environment and, $B_1(s)$, $B_2(s)$ and $B_3(s)$ are standard one-dimensional Brownian motions \citep{rao2014}. It is important to note that the system dynamics (\ref{1}) is a very general case of a standard SIR model. $S^*$, $I^*$ and $R^*$ represent the steady state levels of the state variables in this system.

\begin{as}\label{a0}
	The following set of assumptions regarding the objective function is considered:
	\begin{itemize}
		\item $\{\mathcal F_s\}$ takes the values from a set $\mathcal X\subset \mathbb R^4$. $\{\mathcal F_s\}_{s=0}^t$ is an exogenous Markovian stochastic processes defined on the probability space $(\mathcal X_\infty,\mathcal F_0,\mathcal P)$ where, $\mathcal P$ is the probability measure and $\mathcal X_\infty$ is the functional state space where each function is coming from a smooth manifold.
		
		\item For all $\{e(s),v(s),\be(s),S(s),I(s),R(s)\}$, there exists an optimal vaccination rate and lock-down intensity $\{e^*(s),v^*(s)\}_{s=0}^t$, with initial conditions $\be(0),S(0),I(0)$ and $,R(0)$, which satisfy the stochastic dynamics represented by the equations (\ref{0}) and (\ref{1}) for all continuous time $s\in[0,t]$. 
		
		\item The function $\exp(-rs)\left[S(s)\left(\frac{1}{2}\a_{11}v^2(s)+\a_{12}v(s)+\a_{13}\right)+I(s)\left(\frac{1}{2}\a_{21}e^2(s)+\a_{22}e(s)+\a_{23}\right)\right]\\+\be(e(s),v(s))S(s)I(s)$ is uniformly bounded, continuous on both the state and control spaces and, for a given $\{e(s),v(s),\be(s),S(s),I(s),R(s)\}$, they are $\mathcal P$-measurable.
		
		\item The function $\exp(-rs)\left[S(s)\left(\mbox{$\frac{1}{2}$}\a_{11}v^2(s)+\a_{12}v(s)+\a_{13}\right)+I(s)\left(\mbox{$\frac{1}{2}$}\a_{21}e^2(s)+\a_{22}e(s)+\a_{23}\right)\right]\\+\be(e(s),v(s))S(s)I(s)$ is strictly convex with respect to the state and the control variables.
		
		\item There exists an $\varepsilon>0$ such that for all $\{e(s),v(s),S(s),I(s),R(s)\}$,
		\begin{multline*}
		\E_0\left\{\biggr[\exp(-rs)\left[S(s)\left(\mbox{$\frac{1}{2}$}\a_{11}v^2(s)+\a_{12}v(s)+\a_{13}\right)\right.\right.\\\left.\left.+I(s)\left(\mbox{$\frac{1}{2}$}\a_{21}e^2(s)+\a_{22}e(s)+\a_{23}\right)\right]+\be(e(s),v(s))S(s)I(s)\biggr]\biggr|\mathcal F_0\right\}\geq\varepsilon.
		\end{multline*}
	\end{itemize}
\end{as}

Above assumption guarantees the integrability of the cost function.

\begin{definition}\label{d0}
	For an agent, the optimal state variables $e^*(s),v^*(s),S^*(s),I^*(s)$ and, $R^*(s)$  and their continuous optimal lock-down intensity  $e^*(s)$ and vaccination rate $v^*(s)$ constitute a stochastic dynamic  equilibrium  such that for all $s\in[0,t]$ the conditional expectation of the cost function is	
	\begin{multline*}
	\E_0\left\{\biggr[\exp(-rs)\left[S^*(s)\left(\mbox{$\frac{1}{2}$}\a_{11}v^{*2}(s)+\a_{12}v^*(s)+\a_{13}\right)\right.\right.\\\left.\left.+I^*(s)\left(\mbox{$\frac{1}{2}$}\a_{21}e^{*2}(s)+\a_{22}e^*(s)+\a_{23}\right)\right]+\be(e^*(s),v^*(s))S^*(s)I^*(s)\biggr]\biggr|\mathcal F_0^*\right\}\\
	\leq \E_0\left\{\biggr[\exp(-rs)\left[S(s)\left(\mbox{$\frac{1}{2}$}\a_{11}v^2(s)+\a_{12}v(s)+\a_{13}\right)\right.\right.\\\left.\left.+I(s)\left(\mbox{$\frac{1}{2}$}\a_{21}e^2(s)+\a_{22}e(s)+\a_{23}\right)\right]+\be(e(s),v(s))S(s)I(s)\biggr]\biggr|\mathcal F_0\right\},
	\end{multline*}
	with the dynamics explained in Equations (\ref{0}) and (\ref{1}), where $\mathcal F_0^*$ is the optimal filtration starting at time $0$ so that, $\mathcal F_0^*\subset\mathcal F_0$.
\end{definition}

Define $\mathbf X(s)=[\be(s),S(s),I(s),R(s)]^T$ where $T$ represents the transposition of a matrix such that the dynamic cost function is 
\begin{multline*}
c[u(s),\mathbf X(s)]=\exp(-rs)\left[S(s)\left(\mbox{$\frac{1}{2}$}\a_{11}v^2(s)+\a_{12}v(s)+\a_{13}\right)\right.\\\left.+I(s)\left(\mbox{$\frac{1}{2}$}\a_{21}e^2(s)+\a_{22}e(s)+\a_{23}\right)\right]+\be(e(s),v(s))S(s)I(s),
\end{multline*}
where $u(s)=[e(s),v(s)]^T$. Furthermore, for continuous time $s\in[0,t]$ define
\[
M[\mathbf X(s)]=\inf_{\mathbf X\in\mathcal X}\E_0\left\{\int_0^t c[u(s),\mathbf X(s)]ds\bigg|\mathcal F_0\right\},
\]
where $\mathcal X$ is assumed to be a convex set of state variables.

In the following Proposition \ref{p0}, we will prove the existence of a solution for dynamic cost minimization under complete and perfect information.

\begin{prop}\label{p0}
	Let $c$ be a dynamic quadratic cost function satisfying Assumption \ref{a0} and 
	\[
	\liminf_{\mathbf x\ra\infty}\frac{M[\mathbf X(s)]}{\mathbf X(s)}\geq 0.
	\] 
	Then under perfect and complete information about the pandemic and for all $\mathbf X(s)>0$, there exists a unique solution $\mathbf X^*$ to the problem described in Equation (\ref{0}).
\end{prop}

\begin{proof}
	See the Appendix.
\end{proof}

\begin{remark}
	The condition $\liminf_{\mathbf X\ra\infty}M[\mathbf X(s)]/\mathbf X(s)\geq 0$ in the Proposition \ref{p0} looks strange at first but from the proof we know that it is a necessary condition for existence and uniqueness of the solution of Equation (\ref{0}).
\end{remark}

\subsection{Spread of the Pandemic}

\medskip

\begin{figure}[H]
	\centering
	\includegraphics[width=.5\linewidth]{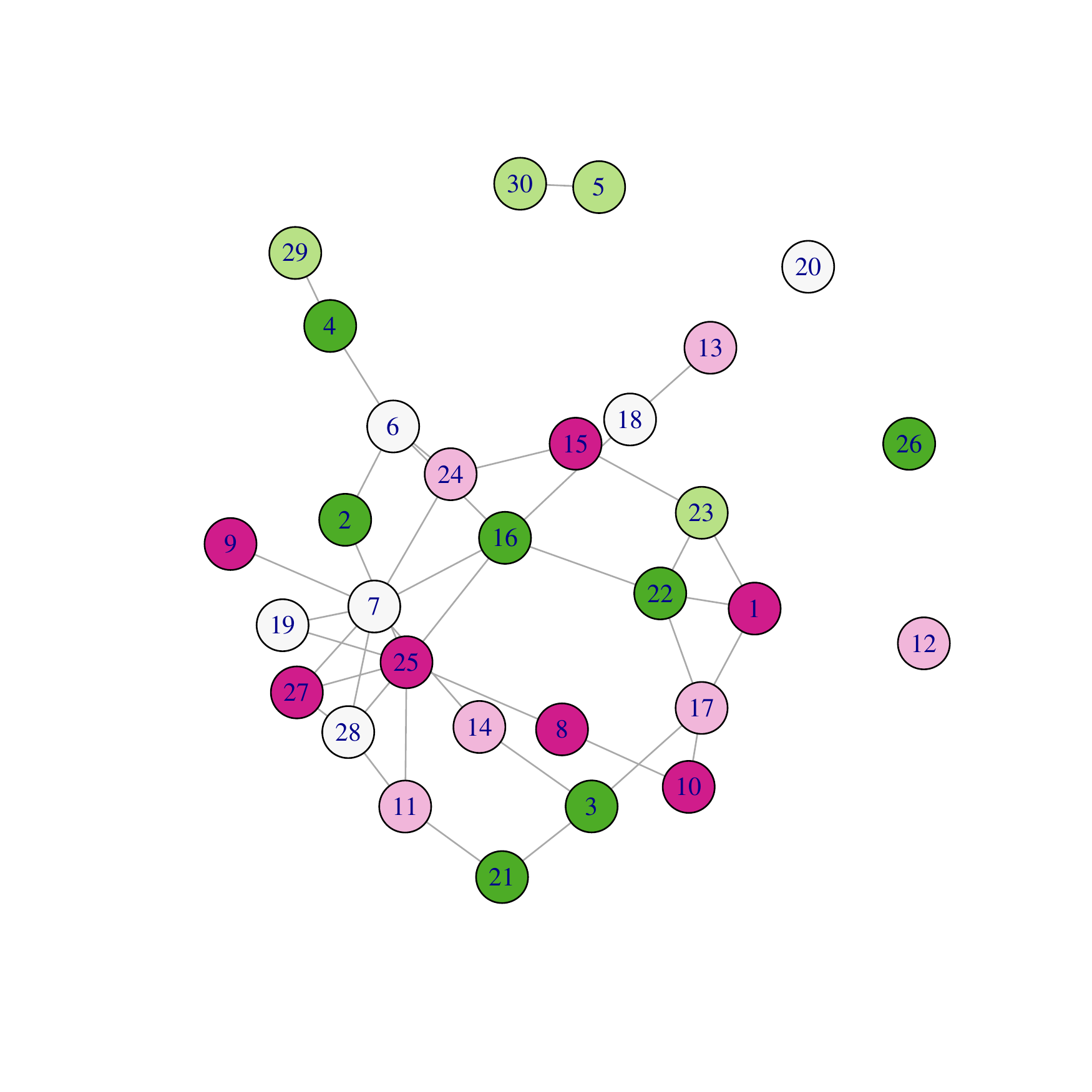}
	\caption{An Erdos-Renyi random network of $30$ agents with five different immunity levels.}
	\label{fig0}
\end{figure}

\begin{figure}[H]
	\begin{minipage}[b]{0.5\linewidth}
		\centering
		\includegraphics[width=\textwidth]{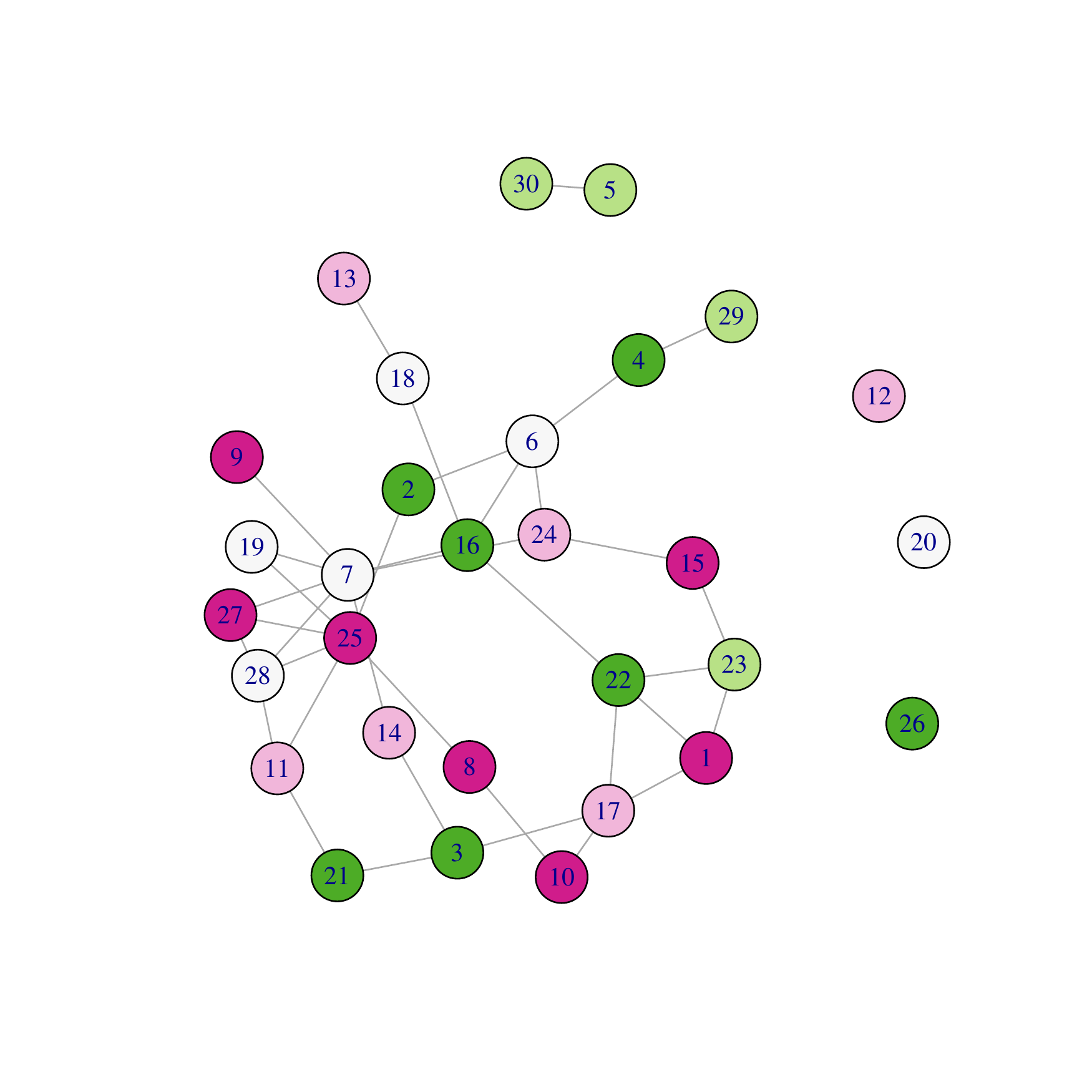}
		\caption{Tie between agents $16$ and $25$ is removed randomly.}
		\label{fig1}
	\end{minipage}
	\hspace{0.5cm}
	\begin{minipage}[b]{0.5\linewidth}
		\centering
		\includegraphics[width=\textwidth]{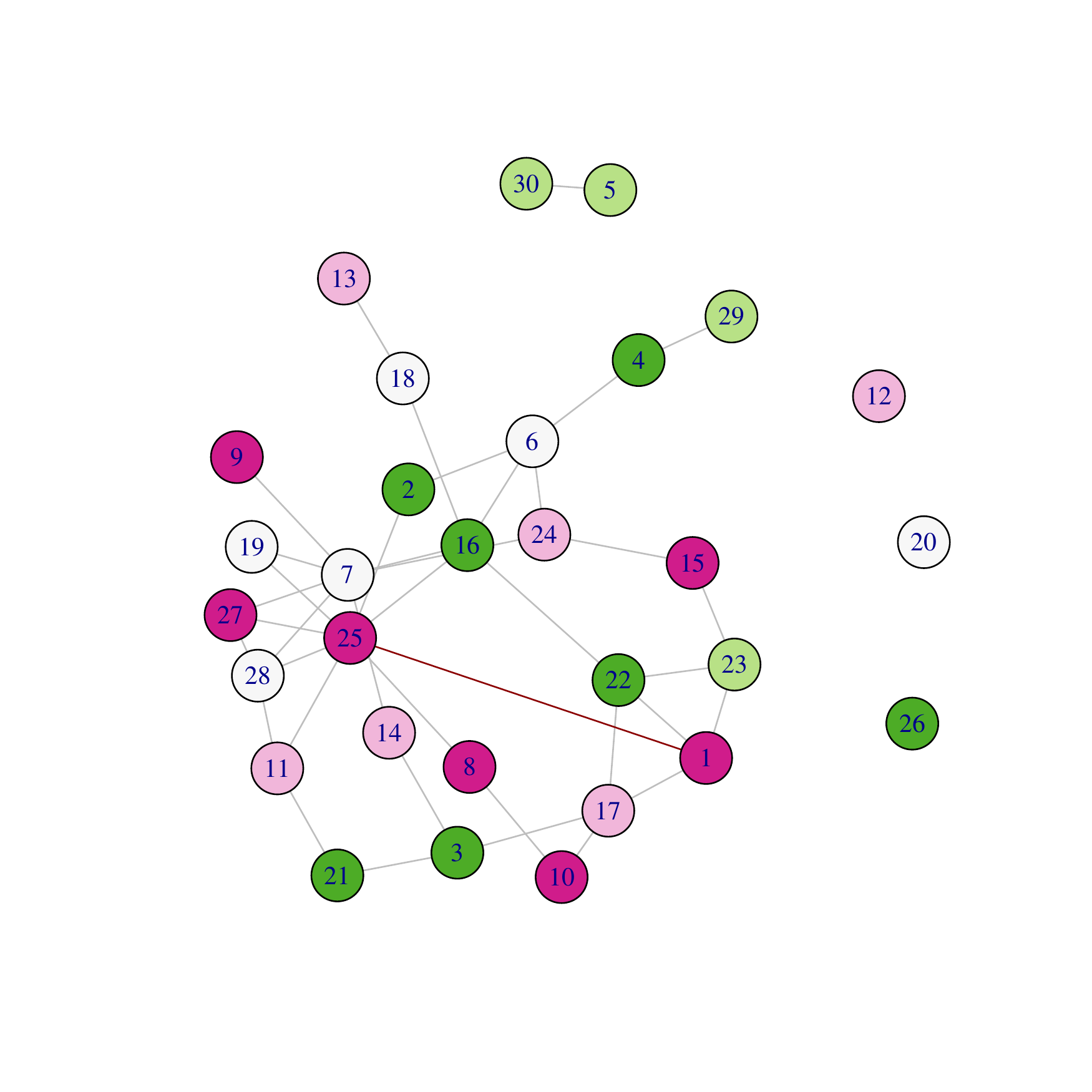}
		\caption{Tie between agents $1$ and $25$ is added randomly.}
		\label{fig2}
	\end{minipage}
\end{figure}

In this section, we are going to discuss the spread of COVID-19 due to social interactions and different levels of immunity levels among humans. We know the immune system is the best defense because it supports the body’s natural ability to defend against pathogens such as viruses, bacteria, fungi, protozoan, and worms, and resists infections \citep{chowdhury2020}. As long the immunity level of a human is properly functional, infections due to a pandemic like COVID-19 go unnoticed. There are three main types of immunity levels such as innate immunity (rapid response), adaptive immunity(slow response), and passive immunity \citep{chowdhury2020}. To determine the interaction among people with different levels of immunity we randomly consider a network of $30$ people. Furthermore, we characterize the immunity levels among five categories as \emph{very low, somewhat low, medium, somewhat high} and \emph{very high}. The subcategories \emph{somewhat high} and \emph{very high} go under innate immunity and, subcategories \emph{very low} and \emph{somewhat low} go under adaptive immunity. We keep passive immunity as \emph{medium} category. We did not subdivide this category under two types: natural immunity, received from the maternal side, and artificial immunity, received from medicine \citep{chowdhury2020} as it is beyond the scope of this paper. In Figure \ref{fig0} we have created an Erdos-Renyi random network of $30$ agents where deep magenta, light magenta, white, lighter green, and deep green represent an agent with \emph{very low, somewhat low, medium, somewhat high} and \emph{very high} respectively.  
In Figure \ref{fig0} let us consider the interaction of agent $25$. According to our setting, this individual has the lowest level of immunity against the pandemic. As the information is perfect and complete everybody in the network has information about COVID-19. Agent $25$ is connected with agents $ 2, 8, 11, 16, 19, 27$ and $28$ where agents $2$ and $16$ have the highest level of immunity system. Assume COVID-19 hits this network and agent $25$ got infected. This person is going to be isolated from some of his adjacent ties, based on a probability-weighted by the level of dissimilarities among their immune systems. Furthermore, agent $25$ would stay with some other non-adjacent agents. In Figure \ref{fig1} we randomly remove the tie between agents $25$ and $16$ who are perfectly opposite in terms of their immune systems. On the other hand, in Figure \ref{fig2} we randomly add a new tie of agent $25$ with a previously nonadjacent agent $1$. Intuitively, one can think about because of COVID-19, agents with similar immune systems tend to come closer.

\begin{figure}[H]
	\centering
	\includegraphics[width=.6\linewidth]{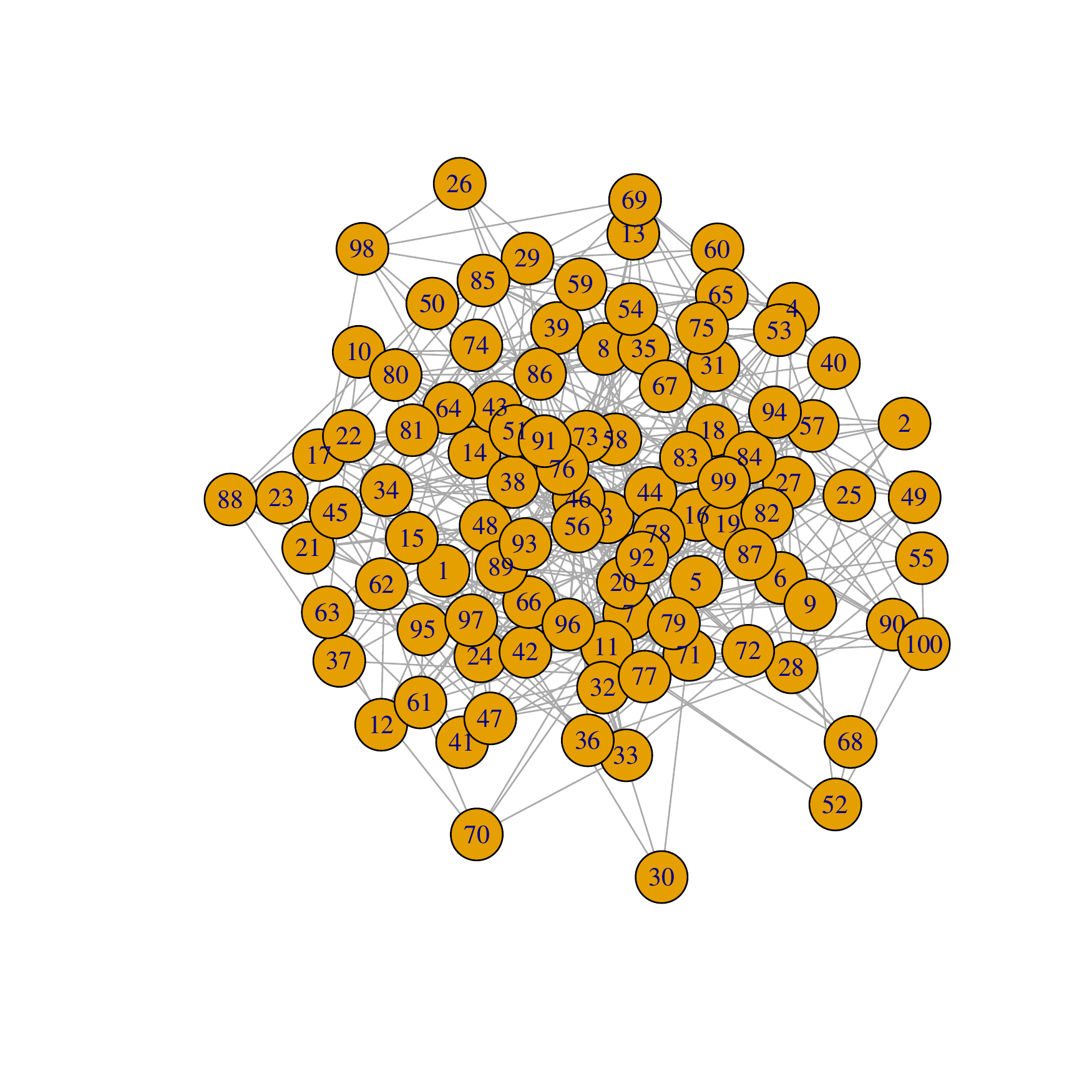}
	\caption{An Erdos-Renyi random network of $100$ agents with five different immunity levels.}
	\label{fig3}
\end{figure}

The temperature takes an important role in spreading the pandemic. If the temperature is high, more people tend to come outside the home and interact. As a result, the spread of the disease would be faster. In order to see interactions between agents in a large network we choose an Erdos-Renyi random network with $100$ agents where $21$ agents have \emph{very low}, $24$ have \emph{somewhat low}, $18$ have \emph{medium}, $20$ have \emph{somewhat high} and $17$ have \emph{very high} immunity systems respectively. Figure \ref{fig3} shows this type of network where agents are connected and disconnected randomly over time based on probabilities weighted by dissimilar immunity levels and temperature of that region.

\begin{figure}[H]
	\centering
	\includegraphics[width=.8\linewidth]{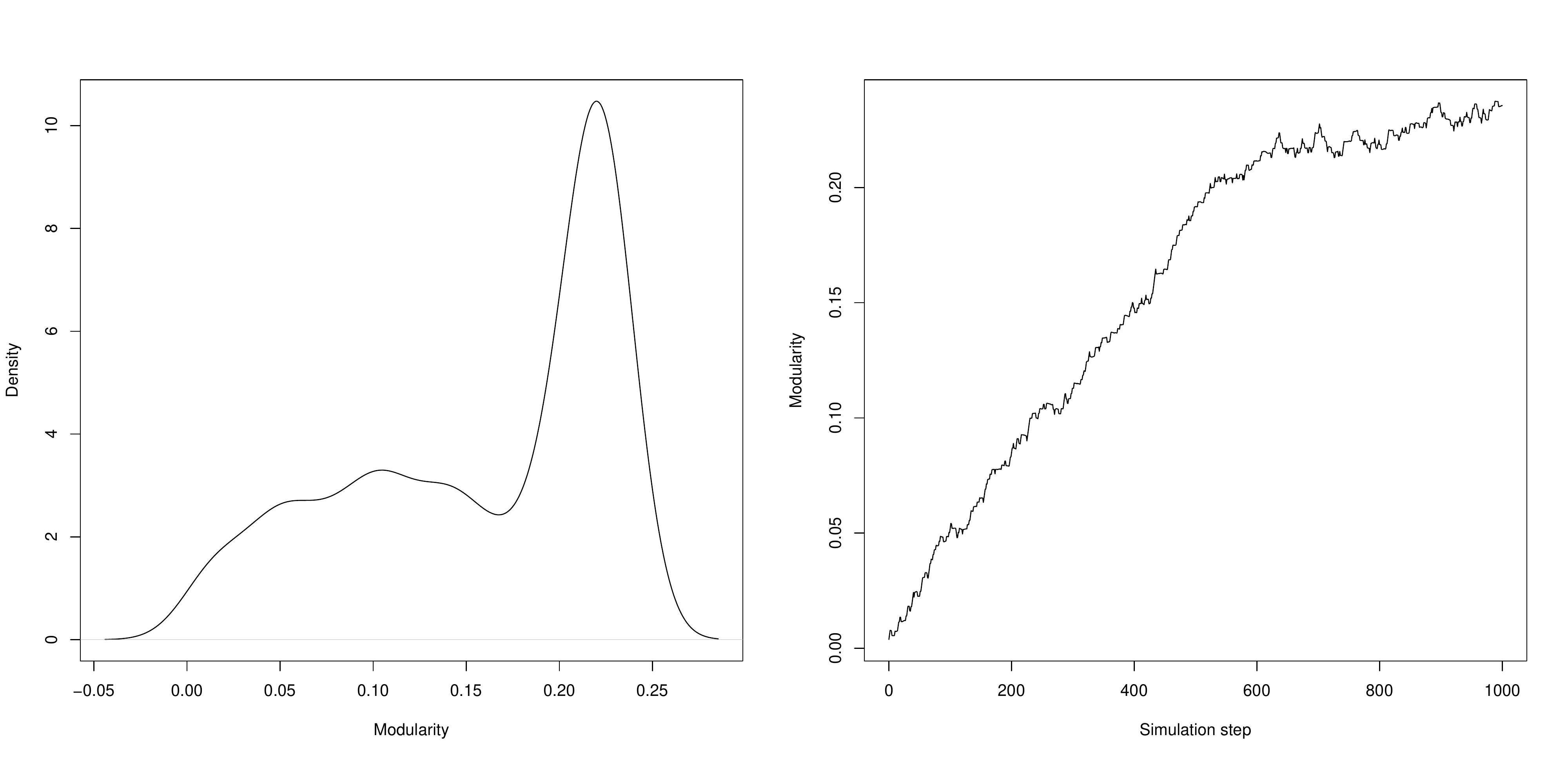}
	\caption{An Erdos-Renyi random network of $100$ agents with $1000$ updates.}
	\label{fig4}
\end{figure}

To create Figure \ref{fig3} an abstract notion of time is used. Here edges are selected randomly for updating assuming that time has some passage between each update. Firstly, the updating function starts with a list of objects that would be used to store an updated network \citep{luke2015, polansky2021}. Then inside a loop a random node is selected, and the update function is called when an existing edge is removed, and a new edge is added. This procedure has two limitations. First, the loop can be replaced by a vectorized function and in each step, this update function stores the entire network which results in very large objects being returned \citep{hua2019,luke2015,pr4}. In Figure \ref{fig4} we update this large network $1000$ times. Before starting the updating process we assume this random network would be more homophilous over time as the updates of edges are partially driven by the similarity of the immunity levels between two agents.

The right panel of Figure \ref{fig4} shows that modularity is lower at the starting network compared to the final network with $1000$ updates whereas the left panel shows the density of this infection network. Intuitively, one might think about after the first incidence of the pandemic a greater part of the network is segregated. Higher modularity at the end implies edges between vertices of similar immunity levels are more likely than edges between different immunity levels. 

Based on the above discussions we are going to construct a stochastic differential equation of the transmission rate of the pandemic $\be$. Consider an Erdos-Renyi random network with the total number of vertices $V$ and edges $\mathcal E$ such that the graph is denoted as $G(V,\mathcal E)$. Let $A(s)$ be the adjacency matrix with each element $a_{m_1m_2}$ for agent $m_1$ and $m_2$. We define the modularity of this network as
\begin{equation*}
\mathcal Q:=\frac{1}{2\mathcal E}\sum_{m_1,m_2}\left[a_{m_1m_2}-\frac{b_{m_1}b_{m_2}}{2\mathcal E}\right]\delta(c_{m_1},c_{m_2}),
\end{equation*}
where $b_{m_i}$ is the degree of the vertex $m_i$ (i.e. agent $m_i$) for all $i=1,2$, $c_{m_i}$ is the community corresponding to $m_i$ with Kronecker delta function $\delta(.,.)$ such that if two different communities merge, $\delta$ takes the value of $1$. As we know, higher temperature increases the transmission rate and, higher lock-down intensity and vaccination rate reduce the transmission rate, the stochastic differential equation is
\begin{equation}\label{3}
d\be(s)=\mathcal QI(s)\left[\be_0T(s)+\be_1M[1-e(s)]^{\theta_1}-\be_2v^{\theta_2}(s)\right]ds+\sigma_4[\be(s)-\be^*]MdB_4(s),
\end{equation}
where $\be\in(0,1)$ are coefficients, $\theta_l>1$ for all $l=1,2$ make the transmission function $\be(v,e)$ a convex function of $e$ and $v$. Moreover, $\be_0$ is the minimum level of infection risk produced if only the essential activities are open, $\be_1$ is the increment in the level of infection, $\be_2$ is the reduction in the level of infection due to vaccination, $M$ is fine particulate matter ($PM_{2.5}>12 \mu g/m^3$) which is an air pollutant and have significant contribution to degrade a person's health,  $\sigma_4$ is a known diffusion coefficient infection dynamics and $dB_4(s)$ is a one-dimensional standard Brownian motion of $\be(e,v)$ with steady state at $\be^*$ and $T(s)$ is the temperature in that region at time $s$.

\subsection{Stochastic SIR Dynamics}

\medskip

For a complete  probability space  $(\mathcal X_\infty,\mathcal F_0,\mathbb P)$ with filtration starting from $\{\mathcal F_s\}_{0\leq s\leq t}$, let \\$\mathbf X(s)=\left[\be(s),S(s),I(s),R(S)\right]^T$ with $L^2$-norm $||\mathbf X(s)||_2=\sqrt{\be^2(s)+S^2(s)+I^2(s)+R^2(s)}$. Let \\$C^{2,1}(\mathbb R^4\times(0,\infty),\mathbb R_+)$ be a family of all non-negative functions $\mathbf Z(s,\mathbf X)$ defined on $\mathbb R^4\times (0,\infty)$ so that they are  twice continuously differentiable in $\mathbf X$ and once in $s$. Define a differential operator $\mathfrak D$ associated with 4-dimensional stochastic differential equation explained in the system of Equations (\ref{1}) and (\ref{3}) as
\begin{equation}\label{5}
d\mathbf X(s)=\bm \mu(s,u,\mathbf X)ds+\bm\sigma(s,\mathbf X) d\mathbf B(s),
\end{equation}
so that
\begin{equation*}
\mathfrak D=\frac{\partial}{\partial s}+\sum_{j=1}^4\bm\mu_{j}(s,u,\mathbf X)\frac{\partial}{\partial X_{_j}}+\frac{1}{2}\sum_{j=1}^4\sum_{j'=1}^4\left[\left[\bm\sigma^T(s,\mathbf X)\bm\sigma(s,\mathbf X)\right]_{jj'}\frac{\partial^2}{\partial X_{j}X_{j'}}\right],
\end{equation*}
where
\[\bm\mu=\begin{bmatrix}
\eta N(s)-\be(e,v)\frac{S(s)I(s)}{\left[1+\rho I(s)\right]+\eta N(s)}-\kappa S(s)-v(s)+\zeta R(s)\\
\be(e,v)\frac{S(s)I(s)}{\left[1+\rho I(s)\right]+\eta N(s)}-(\mu+\kappa)I(s)-e(s)\\
\mu v(s)I(s)-[\kappa +\zeta]e(s) R(s)\\
\mathcal QI(s)\left[\be_0T(s)+\be_1M[1-e(s)]^{\theta_1}-\be_2v^{\theta_2}(s)\right]
\end{bmatrix}\]
and,
\[\bm\sigma=\begin{bmatrix}
\sigma_1(S-S^*) & 0 &0 &0\\
0& \sigma_2(I-I^*) & 0&0\\
0& 0& \sigma_3(R-R^*)& 0\\
0& 0& 0&\sigma_4(\be-\be^*)M
\end{bmatrix}.\]
If the differential operator $\mathfrak D$ operates on a function $\mathcal Z\in C^{2,1}(\mathbb R^4\times(0,\infty);\mathbb R_+)$, such that
\begin{equation}\label{4}
\mathfrak{D}\mathcal Z(s,\mathbf X)=\frac{\partial}{\partial s}\mathcal Z(s,\mathbf X)	+\bm\mu(s,u,\mathbf X)\frac{\partial}{\partial\mathbf X}\mathcal Z(s,\mathbf X)+\frac{1}{2}\text{$trace$}\left\{\bm\sigma^T(s,\mathbf X)\left[\frac{\partial^2}{\partial X^T \partial X} \mathcal Z(s,\mathbf X)\right]\bm\sigma(s,\mathbf X)\right\},	
\end{equation}	
where T represents a transposition of a matrix.
\begin{as}\label{a1}
	For $t>0$, let $\bm{\mu}(s,u,\mathbf X):[0,t]\times [0,1]^2\times \mathbb{R}^4 \ra\mathbb{R}$ and ${\bm\sigma(\mathbf X)}: \mathbb{R}^4 \ra\mathbb{R}$ be some measurable function and, for some positive constant $K_0$,  $\mathbf X\in\mathbb{R}^4$ we have linear growth as
	\[
	|\bm{{\mu}}(s,u,\mathbf X)|+
	|{\bm\sigma(\mathbf X)}|\leq 
	K_0(1+|\mathbf X|),
	\]
	such that, there exists another positive, finite, constant $K_1$ and for a different state variable vector ${\widetilde {\mathbf X}}$ 
	such that the Lipschitz condition,
	\[
	|\bm{{\mu}}(s,u,\mathbf X)-\bm{{\mu}}(s,u,\widetilde{\mathbf X})|+|\bm\sigma(\mathbf X)-\bm\sigma_0^k(\widetilde{\mathbf X})|\leq K_1\ |\mathbf X-{\widetilde {\mathbf X}}|,\notag
	\]
	$ {\widetilde {\mathbf X}}\in\mathbb{R}^4$ is satisfied and
	\[
	|\bm{{\mu}}(s,u,\mathbf X)|^2+
	|\bm\sigma(\mathbf X)|^2\leq K_1^2
	(1+|{\widetilde {\mathbf X}}|^2).
	\]
\end{as}

\begin{as}\label{a2}
	Assume $(\mathcal X_\infty,\mathcal F_0,\mathcal P)$ is a stochastic basis where the filtration $\{\mathcal F_s\}_{0\leq s\leq t}$ supports a 4-dimensional Brownian motion $\mathbf B(s)=\{\mathbf B(s)\}_{0\leq s\leq t}$. $\mathcal F_0$ is the collection of all $\mathbb R$-values progressively measurable process on $[0,t]\times\mathbb R^4$ and the subspaces are 
	\[
	\mathbb F^2:=\left\{ \mathbf X\in\mathcal F_0;\ \E_0\int_0^t|\mathbf X(s)|^2ds<\infty\right\}
	\] 
	and,
	\[
	\mathbb S^2:=\left\{\mathbf X\in\mathcal F_0;\ \E_0\sup_{0\leq s\leq t}|\mathbf X(s)|^2<\infty\right\},
	\]
	where $\mathcal X_\infty$ is a Borel $\sigma$-algebra and $\mathcal P$ is a probability measure \citep{carmona2016}. Furthermore, the 4-dimensional Brownian motion corresponding to the vector of state variables in this system is defined as 
	\[
	\mathbf B:=\left\{\mathbf X\in\mathcal F_0;\ \sup_{0\leq s\leq t}|\mathbf X(s)|<\infty;\ \mathcal P-a.s.\right\}.\]
\end{as}

\begin{prop}\label{p1}
	Consider a small continuous time interval $[s,\tau]\in(0,t)$. Also assume the left hand side of the Equation (\ref{4}) is zero and $\mathcal Z$ solves the Cauchy problem (\ref{4}) with terminal condition $\mathcal Z(\tau,\mathbf X)=\Phi(\mathbf X)$. Let pandemic state variables $\mathbf X$ follows the stochastic differential equation (\ref{5}) such that 
	\[
	\bm\sigma(s,\mathbf X)\frac{\partial}{\partial \mathbf X}\mathcal Z(s,\mathbf X)\in L^2,\ \text{for all}\ s\leq\tau,\ \mathbf X\in\mathbb R^4,
	\]
	then $\mathcal Z$ has the stochastic Feynman-Kac representation
	\[
	\mathcal Z(s,\mathbf X)=\E_s[\Phi(\mathbf X(\tau))].
	\]
\end{prop}

\begin{proof}
	See the Appendix.
\end{proof}

\begin{prop}\label{p2}
	If Assumptions \ref{a1} holds then under complete information about the pandemic and for continuous $s\in[0,t]$ the pandemic dynamics expressed in Equation (\ref{5}) has a strong unique solution.
\end{prop}

\begin{proof}
	See the Appendix.
\end{proof}

\begin{prop}\label{p3}
	Let the initial state variable of SIR model with stochastic infection  $\mathbf X(0)\in\bf L^2$	is independent of Brownian motion $\mathbf B(s)$ and the drift and diffusion coefficients $\bm\mu(s,u,\mathbf X)$ and $\bm\sigma(s,\mathbf X)$ respectively follow Assumptions \ref{a1} and \ref{a2}. Then the pandemic dynamics in Equation (\ref{5}) is in space of the real valued process with filtration $\{\mathcal F_s\}_{0\leq s\leq t}$ and this space is denoted by $\mathcal{F}_0$. Moreover, for some finite constant $c_0>0$, continuous time $s\in[0,t]$ and Lipschitz constants $\bm\mu$ and $\bm\sigma$, the solution satisfies,
	\begin{equation}\label{6}
	\E \sup_{0\leq s\leq t}|\mathbf X(s)|^2\leq c_0(1+\E|\mathbf X(0)|^2)\exp{(c_0t)}.
	\end{equation}
\end{prop}

\begin{proof}
	See the Appendix.
\end{proof}

Propositions \ref{p1}-\ref{p3} tell us about the uniqueness and measurability of the system of stochastic SIR dynamics with infection dynamics. It is important to know that we assume the information available regarding the pandemic is complete and perfect and that all the agents in the system are risk-averse. Therefore, once a person in a community gets infected by COVID-19, everybody gets information immediately and that agent becomes isolated from the rest. 

\subsection{Main Results}

\medskip

An agent's objective is to minimize the quadratic cost function expressed in Equation (\ref{0})  subject to the dynamic system represented by the equations  (\ref{1}) and (\ref{3}). Following \cite{pramanik2020optimization} the quantum Lagrangian of an agent in this pandemic environment is 
\begin{equation}\label{9}
\mathcal L(s,u,\mathbf X):=\E_s\{c\left[u(s),\mathbf X(s)\right]+\lambda\left[\bm \mu(s,u,\mathbf X)ds+\bm\sigma(s,\mathbf X) d\mathbf B(s)-\Delta\mathbf X\right]\},
\end{equation}
where $\Delta \mathbf X=\mathbf X(s+\varepsilon)-\mathbf X(s)$ for all $\nu\in[s,s+\varepsilon]$, $\varepsilon\downarrow 0$ and $\E_s[.]:=\E[.|\mathbf X(s),\mathcal F_s]$. In Equation (\ref{9}) $\lambda>0$ is a time-independent quantum Lagrangian multiplier. At the time $s$ an agent can predict the severity of the pandemic based on all information available regarding state variables at that time; moreover, throughout interval $[s,s+\varepsilon]$ the agent has the same conditional expectation which ultimately gets rid of the integration.

\begin{prop}\label{p4}
	For any two different immunity groups, if the probability measures of getting affected by the pandemic are $\mathcal P_1$ and $\mathcal P_2$ respectively on $H\in(\mathcal X_\infty,\mathcal P,\mathcal F_0)$ so that the total variation difference between $\mathcal P_1$ and $\mathcal P_2$ is
	\begin{align}
	||\mathcal P_1-\mathcal P_2||_{tv}&=\sup\left\{|\mathcal P_1(\mathcal L)-\mathcal P_2(\mathcal L)|;\ \text{for all}\ \mathcal L\in H\right\}\label{p4.0}\\
	&=1-\sup_{\eta\leq\mathcal P_1,\mathcal P_2}\eta(H)\label{p4.1}\\
	&=1-\inf\sum_{k=1}^K\left[\mathcal P_1(B_k)\wedge\mathcal P_2(B_k)\right]\label{p4.2},
	\end{align}
	where $B_k\subset H$ so that $\bigcap_{k=1}^K B_k=\emptyset$ for all $k\in[1,K]$ and $K\geq 1$.	
\end{prop}

\begin{proof}
	See the Appendix.
\end{proof}

\begin{remark}
	In Proposition \ref{p4} $B_k$ is a set of communities of agents such that no two of them never socially interact. Furthermore, Proposition \ref{p4} tells us that if the same variant of COVID-19 hits a community with two agents differed by their immunities, total variation of infection is the suprimum of two infection probabilities of their quantum Lagrangians.
\end{remark}

\begin{prop}\label{p5}
	Suppose, the domain of the quantum Lagrangian $\mathcal L$ is non-empty, convex and compact denoted as $\widetilde\Xi$ such that $\widetilde\Xi\subset\mathcal U\times\mathcal X\subset \mathbb R^{6}$. As $\mathcal L: \widetilde\Xi\ra\widetilde\Xi$ is continuous, then  there exists a vector of state and control variables $\bar Z^*=[v^*,e^*,\be^*,S^*,I^*,R^*]^T$ in continuous time $s\in[0,t]$ such that $\mathcal L$ has a fixed-point in Brouwer sense, where $T$ denotes the transposition of a matrix.
\end{prop}

\begin{proof}
	See the Appendix.
\end{proof}

\begin{theorem}\label{t0}
	Consider an agent's objective is to minimize $M[\mathbf X(s)]$ subject to the stochastic dynamic system explained in the Equation (\ref{5}) such that the Assumptions \ref{a0}-\ref{a2} and Propositions \ref{p0}-\ref{p5} hold. For a $C^2$-function $\tilde f(s,\bar Z)$ and for all $s\in[0,t]$ there exists a function  $g(s,\mathbf X)\in C^2([0,t]\times\mathbb{R}^{4})$ such that  $\widetilde Y=g(s,\mathbf X)$  with  an It\^o process $\widetilde Y$ optimal ``lock-down" intensity $e^*$ and vaccination rate $v^*$ are the solutions 
	\begin{equation}\label{12}
	-\frac{\partial }{\partial u}\tilde f(s,\bar Z)\Psi_s^{\tau}(\mathbf X)=0,
	\end{equation}
	where $\Psi_s^{\tau}$ is some transition wave function in $\mathbb R^{4}$.
\end{theorem}

\begin{proof}
	See the Appendix.
\end{proof}

\begin{remark}
	Proposition \ref{p5} tells us that this pandemic system has a Brouwer fixed point $\bar Z^*$ and as information is perfect and complete, for a given $g(s,\mathbf X)$, this fixed point is unique. Theorem \ref{t0} helps us to determine those fixed points. Since we are assuming feedback controls, once we obtain a steady state $\bar Z^*$, $u^*$ is automatically achieved.
\end{remark}

\section{Computation}
Theorem \ref{t0} determines the solution of an optimal \emph{lock-down} intensity and \emph{vaccination} rate for a generalized stochastic pandemic system. Consider a function $g(s,\mathbf X)\in C^2([0,t]\times\mathbb R^{4})$ such that \citep{rao2014},
\[
g(s,\mathbf X)=[s\be-1-\ln(\be)]+[sS-1-\ln(S)]+[sI-1-\ln(I)]+[sR-1-\ln(R)],
\]
with $\partial g/\partial s=S+I+R+\be$, $\partial g/\partial x_i=s-1/x_i$, $\partial^2 g_/\partial x_i^2=-1/x_i^2$ and $\partial^2 g/\partial x_i\partial x_j=0$, for all $i\neq j$ where $x_i$ is $i^{th}$ state variable of $\mathbf X$ for all $i=1,...,4$ and $\ln$ stands for natural logarithm. In other words, $x_1=\be, x_2=S,x_3=I$ and $ x_4=R$. Therefore,
\begin{multline*}
\tilde f(s,\bar Z)=\exp(-rs)\left[S\left(\mbox{$\frac{1}{2}$}\a_{11}v^2+\a_{12}v+\a_{13}\right)+I\left(\mbox{$\frac{1}{2}$}\a_{21}e^2+\a_{22}e+\a_{23}\right)\right]+\be SI\\
+[s\be-1-\ln(\be)]+[sS-1-\ln(S)]+[sI-1-\ln(I)]+[sR-1-\ln(R)]+(\be+S+I+R)+\left(s-\frac{1}{\be}\right)\\
\times \mathcal QI\left[\be_0T+\be_1M(1-e)^{\theta_1}-\be_2v^{\theta_2}\right]+\left(s-\frac{1}{S}\right)\biggr\{\eta N-\frac{\be SI}{\left(1+\rho I\right)+\eta N}-\kappa S-v+\zeta R\biggr\}\\
+\left(s-\frac{1}{I}\right)\left\{\frac{\be SI}{\left(1+\rho I\right)+\eta N}-(\mu+\kappa)I-e\right\}+\left(s-\frac{1}{R}\right)\left[\mu vI-(\kappa +\zeta)e R\right]\\
-\mbox{$\frac{1}{2}$}\left\{\sigma_1(S-S^*)\frac{1}{S^2}+\sigma_2(I-I^*)\frac{1}{I^2}+\sigma_3(R-R^*)\frac{1}{R^2} +\sigma_4(\be-\be^*)\frac{1}{\be^2}\right\}.
\end{multline*}
To satisfy Equation (\ref{16}), either $\frac{\partial \tilde f}{\partial u}=0$ or $\Psi_s^{\tau}=0$. As $\Psi_s^{\tau}$ is a wave function, it cannot be zero. Therefore, $\frac{\partial \tilde f}{\partial u}=0$ for all $u=\{e,v\}$. Therefore, for $\theta_1=2$ the lock-down intensity is,
\[
e^*=\frac{A_2+A_3}{A_1+A_2},
\]
where $A_1=\exp(-rs)I\a_{21}$, $A_2=2\mathcal QI\be_1M\left(s-\frac{1}{\be}\right)$ and $A_3=\left(s-\frac{1}{I}\right)+R\left(s-\frac{1}{R}\right)(\kappa+\zeta)-\exp(-rs)I\a_{22}$. On the other hand, for $\theta_2=2$ the vaccination rate is,
\[
v^*=\frac{B_3}{B_1-B_2},
\]
where $B_1=\exp(-rs)S\a_{11}$, $B_2=2\mathcal Q I\be_2\left(s-\frac{1}{\be}\right)$ and $B_3=\left(s-\frac{1}{S}\right)-\mu I\left(s-\frac{1}{R}\right)-\exp(-rs)S\a_{12}$ so that $B_1>B_2$. 

\subsection{Simulation Studies}

\begin{figure}[H]
	\centering
	\includegraphics[width=.5\linewidth]{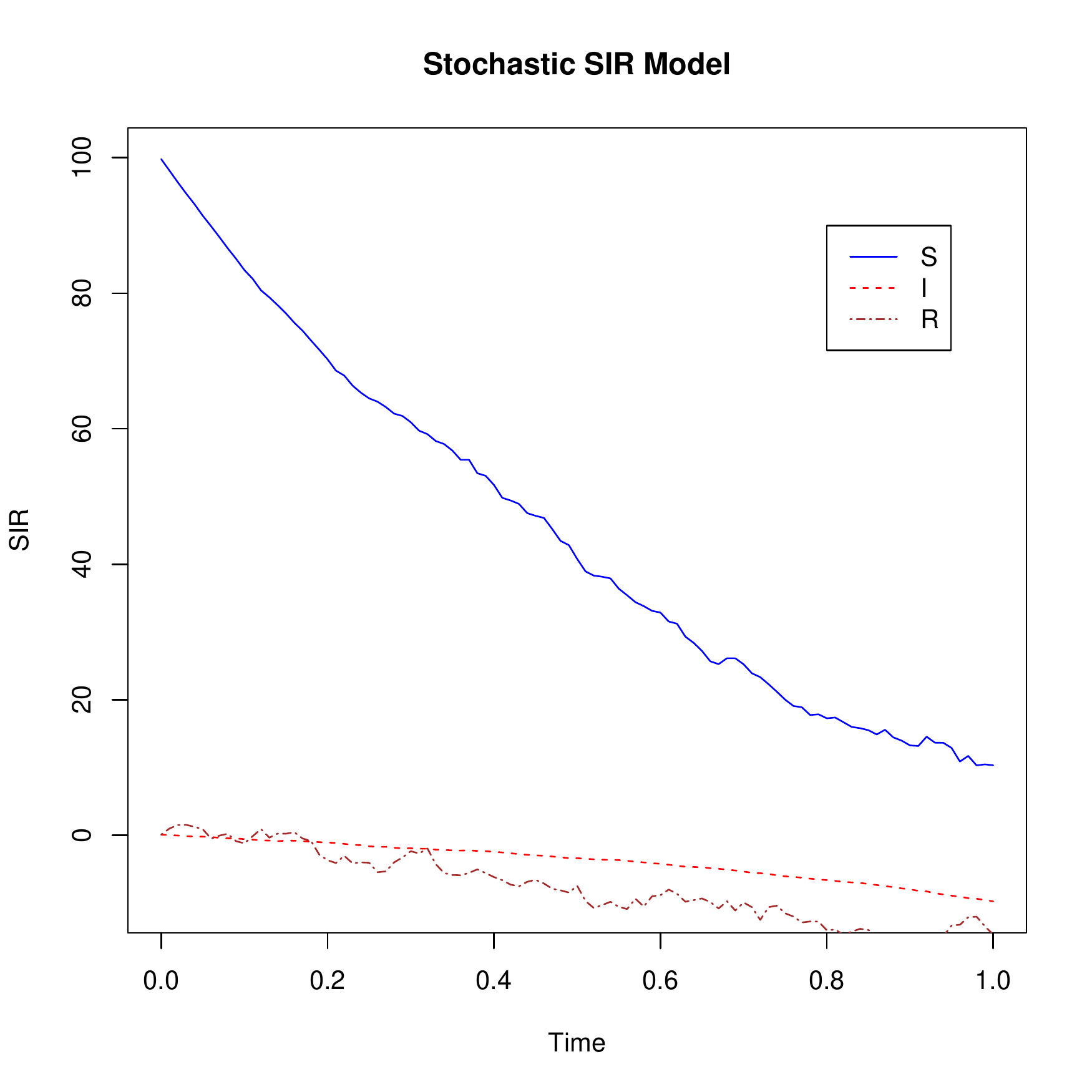}
	\caption{SIR model with higher volatility.}
	\label{fig5}
\end{figure}

Values from Table \ref{tab} have been used to perform simulation studies. These values and initial state variables are obtained from \cite{caulkins2021optimal} and \cite{rao2014}. We did simulate the stochastic SIR model $100$ times with different diffusion coefficients. Figure \ref{fig5} assumes $\sigma_1=0.1$, $\sigma_2=0.06$ and $\sigma_3=0.12$.

\begin{table}[H]
	\centering
	\caption{Parametric values taken from \cite{caulkins2021optimal} and \cite{rao2014}.}   
	
	\begin{tabular}{ |p{3cm}||p{3cm}|p{10cm}| }
		\hline
		\multicolumn{3}{|c|}{Parameter values and initial state variable values.} \\
		\hline
		Variable& Value  &Description\\
		\hline
		$\eta$  & 0.001    & Birth-rate\\
		$\be$& 1& Initial infection\\
		$\be_0$ &   0 & Minimal level of infection   \\
		$\be_1$ & 0.2 & Increment in the level of infection\\
		$\be_3$ & 0.2 & Reduction in the level of infection due to vaccination\\
		$e(0)$    & 1 & Initial lock-down intensity\\
		$\kappa$ &  0.2  & Death-rate\\
		$\zeta$ & 0.001  & Rate by which recovered get susceptible again   \\
		$\mu$ & 0.3  & Natural recovery rate\\
		$\rho$ & 0.5 &Psychological or inhibitory coefficient\\
		$\theta_l$ & 2 & Convexity coefficient of transmission function\\
		$M$ & 12.5 & Fine particulate matter\\
		$\mathcal Q$ & 0.5 & Modularity of network\\
		S(0) & 99.8 &Initial susceptible population\\
		I(0) & 0.1 & Initial infected population\\
		R(0) & 0.1 & Initial recovered population\\
		$v$ & 0.674 & Stable fully vaccination rate\\
		$\a_{ij}$ & $\frac{1}{3}$ & Coefficients of cost function\\
		\hline
	\end{tabular}
	\label{tab}
\end{table}

\begin{figure}[H]
	\begin{minipage}[b]{0.4\linewidth}
		\centering
		\includegraphics[width=\textwidth]{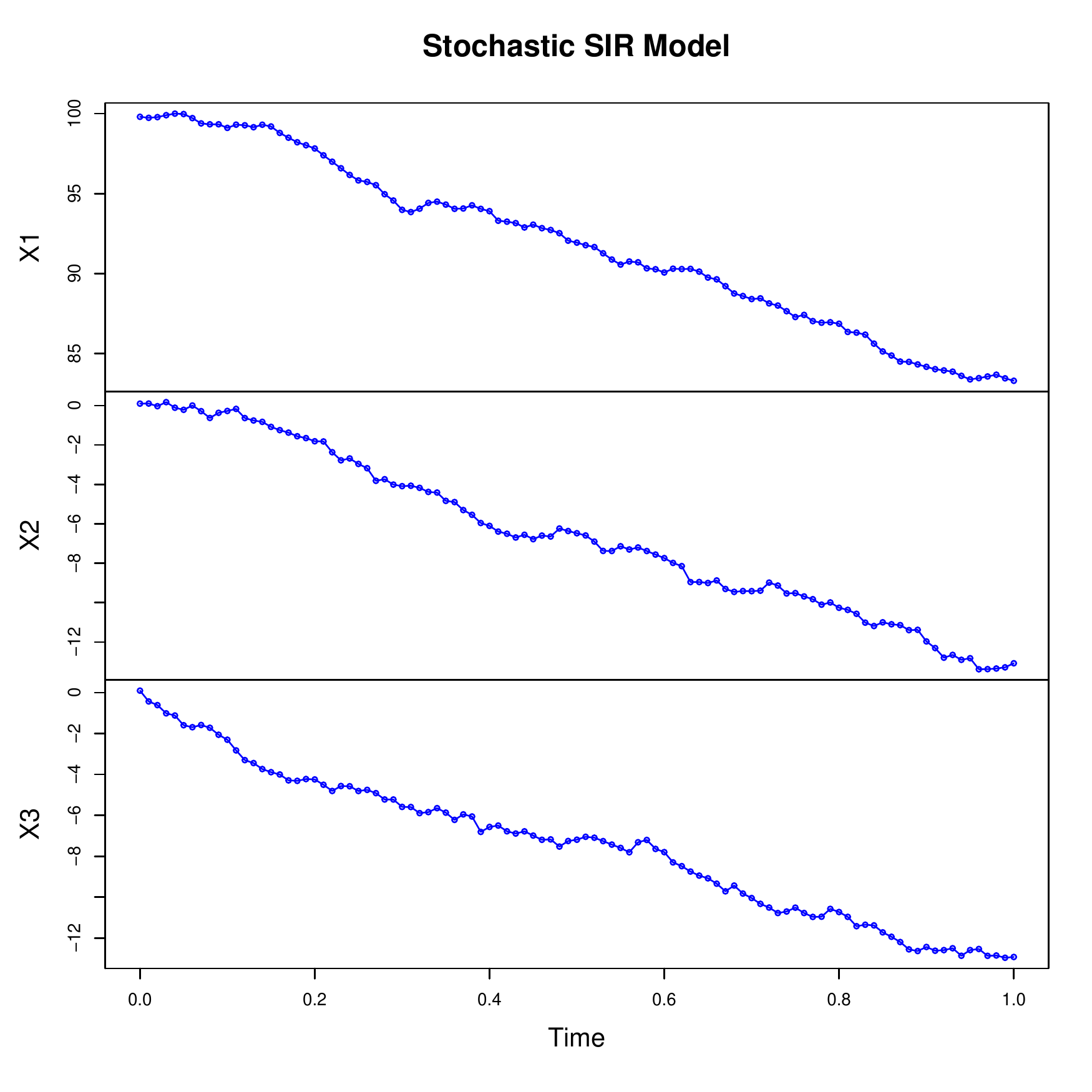}
		\caption{Model with $\sigma_1=0.05$, $\sigma_2=0.01$ and $\sigma_3=0.03$.}
		\label{fig6}
	\end{minipage}
	\hspace{0.5cm}
	\begin{minipage}[b]{0.4\linewidth}
		\centering
		\includegraphics[width=\textwidth]{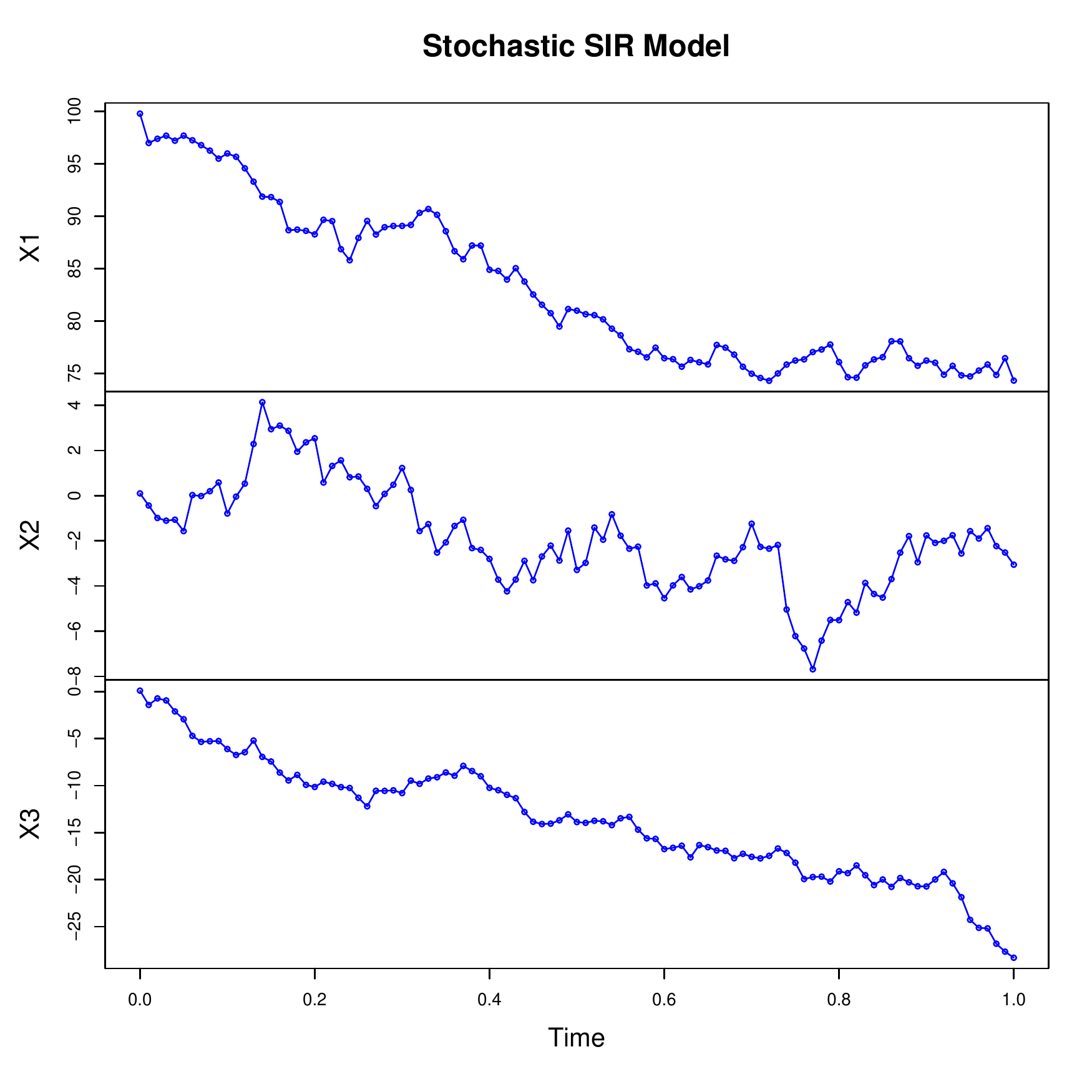}
		\caption{Model with $\sigma_1=0.1$, $\sigma_2=0.06$ and $\sigma_3=0.12$.}
		\label{fig7}
	\end{minipage}
\end{figure}

\begin{figure}[H]
	\begin{minipage}[b]{0.4\linewidth}
		\centering
		\includegraphics[width=\textwidth]{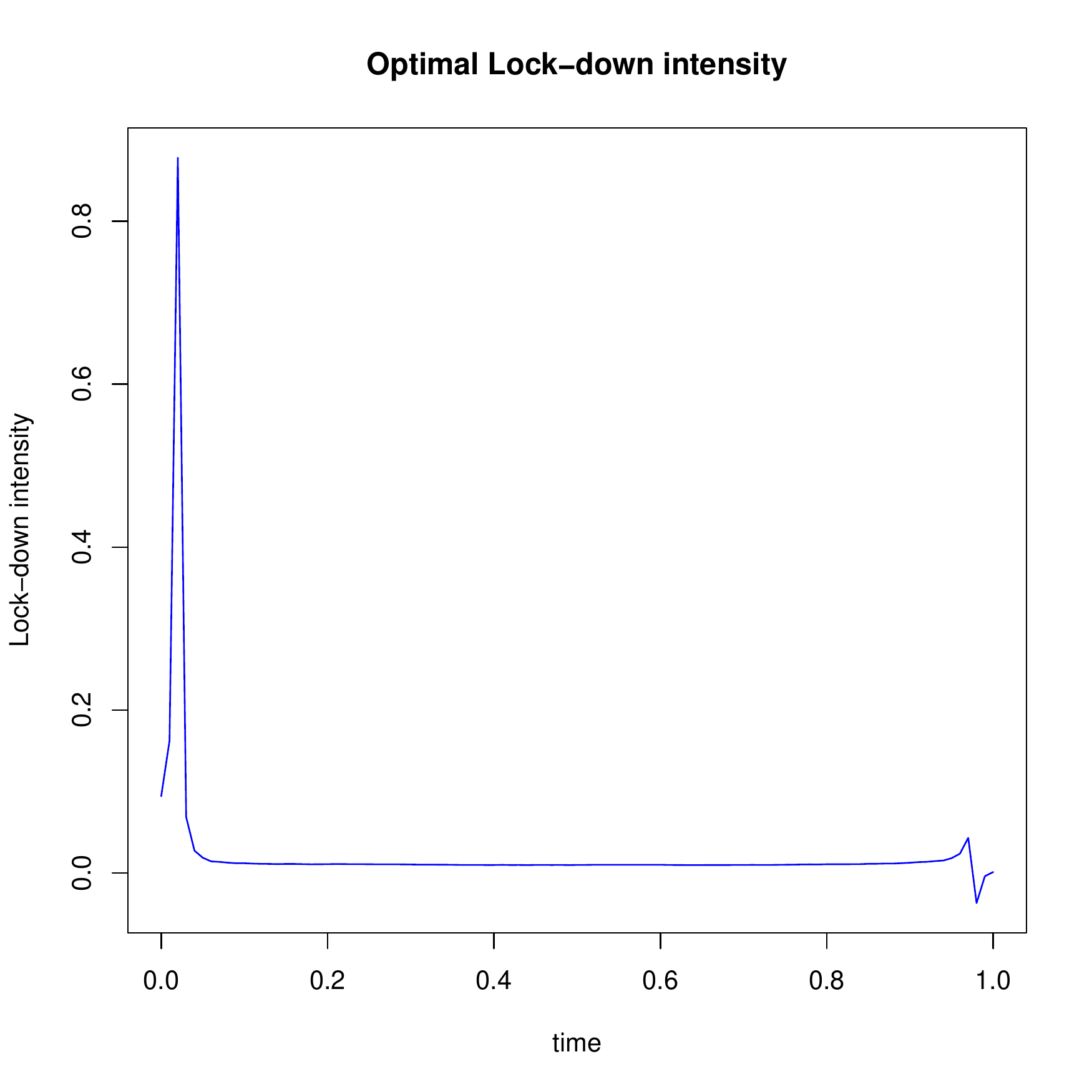}
		\caption{Lock-down with diffusion coefficients $\sigma_2=0.06$ and $\sigma_3=0.12$.}
		\label{fig8}
	\end{minipage}
	\hspace{0.5cm}
	\begin{minipage}[b]{0.4\linewidth}
		\centering
		\includegraphics[width=\textwidth]{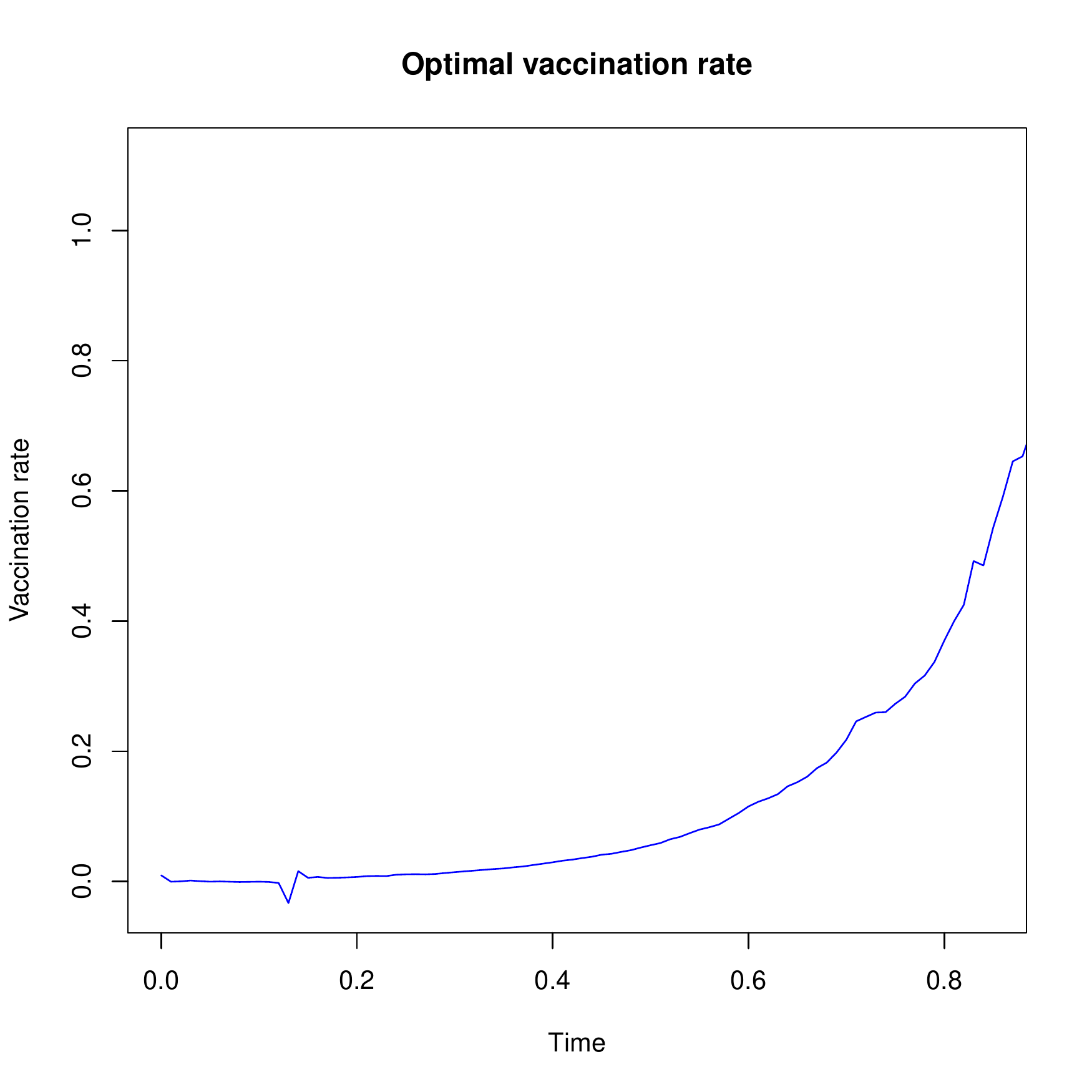}
		\caption{Vaccination with $\sigma_1=0.1$, $\sigma_2=0.06$ and $\sigma_3=0.12$.}
		\label{fig9}
	\end{minipage}
\end{figure}

Since the diffusion coefficients are relatively high, we can see more fluctuations. To observe the behavior of each of the susceptible (S), infected (I), and recovered (R) curves we construct Figures \ref{fig6}-\ref{fig7}. In these figures, $X1$, $X2$, and $X3$ curves represent S, I, and R respectively. When the diffusion coefficients are low then all three curves have a downward pattern as in Figure \ref{fig6}. Once these coefficients increased to $\sigma_1=0.1$, $\sigma_2=0.06$ and $\sigma_3=0.12$, the $X2$ curve in Figure \ref{fig7} starts to behave ergodically, while $X1$ and $X3$ keep their downward trends with more fluctuations. Figures \ref{fig8} and \ref{fig9} represent the behavior of optimal lock-down intensity and vaccination rate over time. Our model says, under higher volatility of the pandemic the vaccination rate is increasing over time because people are \emph{risk-averse} and the information regarding this disease \emph{perfect} and \emph{complete}.

On the other hand, under $\sigma_2=0.06$ and $\sigma_3=0.12$ figure \ref{fig8} implies initially people did not know about the severity of the disease, and therefore, they come outside their homes and work. Slowly they become afraid of being infected and stopped going out and finally, very close to the terminal point the intensity increased because of the high vaccination rate (i.e. Figure \ref{fig9}).

\subsection{Real data analysis}

\medskip

In this section we determine the parametric values from UK data at the beginning of $2021$ \citep{ons2021,steel2021}. The initial conditions of susceptibility, infection and recovery are taken at the beginning of $2021$ (i.e. early January) when a post-Christmas spike of infection took place and the vaccination had just begun \citep{maciejowski2022}. Therefore, $S(0)=84.19$, $I(0)=1.89$ and $R(0)=13.82$. The values of the initial conditions are derived from the \emph{Office for National Statistics} \citep{ons2021} by summing over England, Wales, Scotland, and Northern Ireland and averaging over two months $10$ December 2020 to $6$ January 2021 and, $7$ January to $3$ February 2021. Moreover, initial condition of recovery is calculated by using the formula $R(0)=100-S(0)-I(0)$. The estimate of death rate $\kappa=0.01$ is determined by dividing the cumulative number of deaths up to 14 January 2021 by the estimated number in the recovered category \citep{maciejowski2022}. The birthrate ($\eta$) at this point of time in the UK is $0.0558$, the initial \emph{lock down intensity} $e(0)=0.75$, the vaccination rate with first and second doses are $0.0291$ and $0.00557$ respectively. Since throughout this paper we consider only the full vaccination rate, we are going to use $v=0.00557$. We assume the total population of the UK at that time was $67.22$ million. We also determine increment in the level of infection $\be_1=0.536$ under the assumption of $\be_1=\be_2$, $\sigma_2=0.08557$ and the rate at which recovered agent gets susceptible again $\zeta=0.000152$ \citep{ons2021}. For the other parameter values, we are going to use Table \ref{tab}.

\begin{figure}[H]
	\begin{minipage}[b]{0.4\linewidth}
		\centering
		\includegraphics[width=\textwidth]{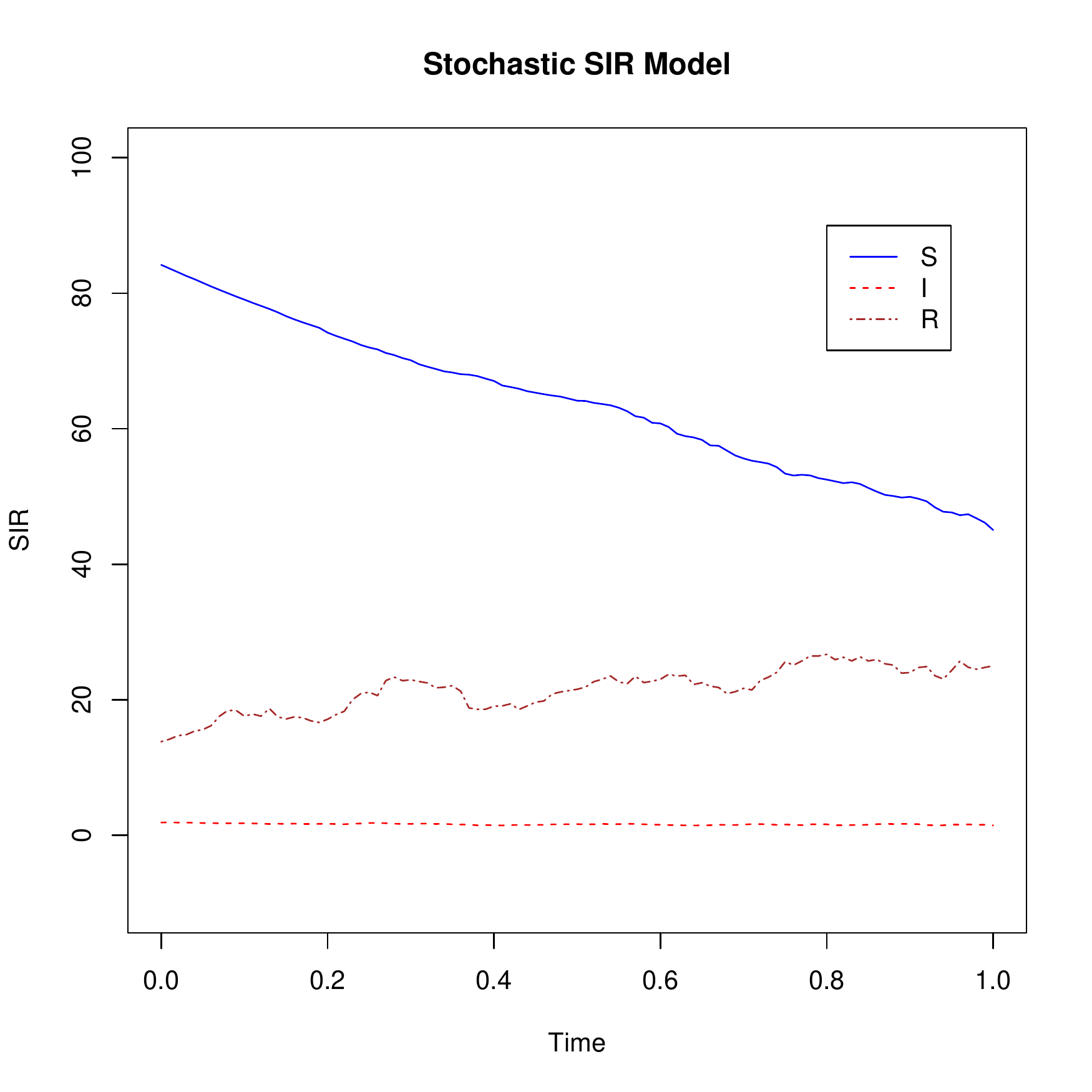}
		\caption{ SIR Model of UK data with $\sigma_1=0.05$, $\sigma_2=0.08557$ and $\sigma_3=0.12$.}
		\label{fig10}
	\end{minipage}
	\hspace{0.5cm}
	\begin{minipage}[b]{0.4\linewidth}
		\centering
		\includegraphics[width=\textwidth]{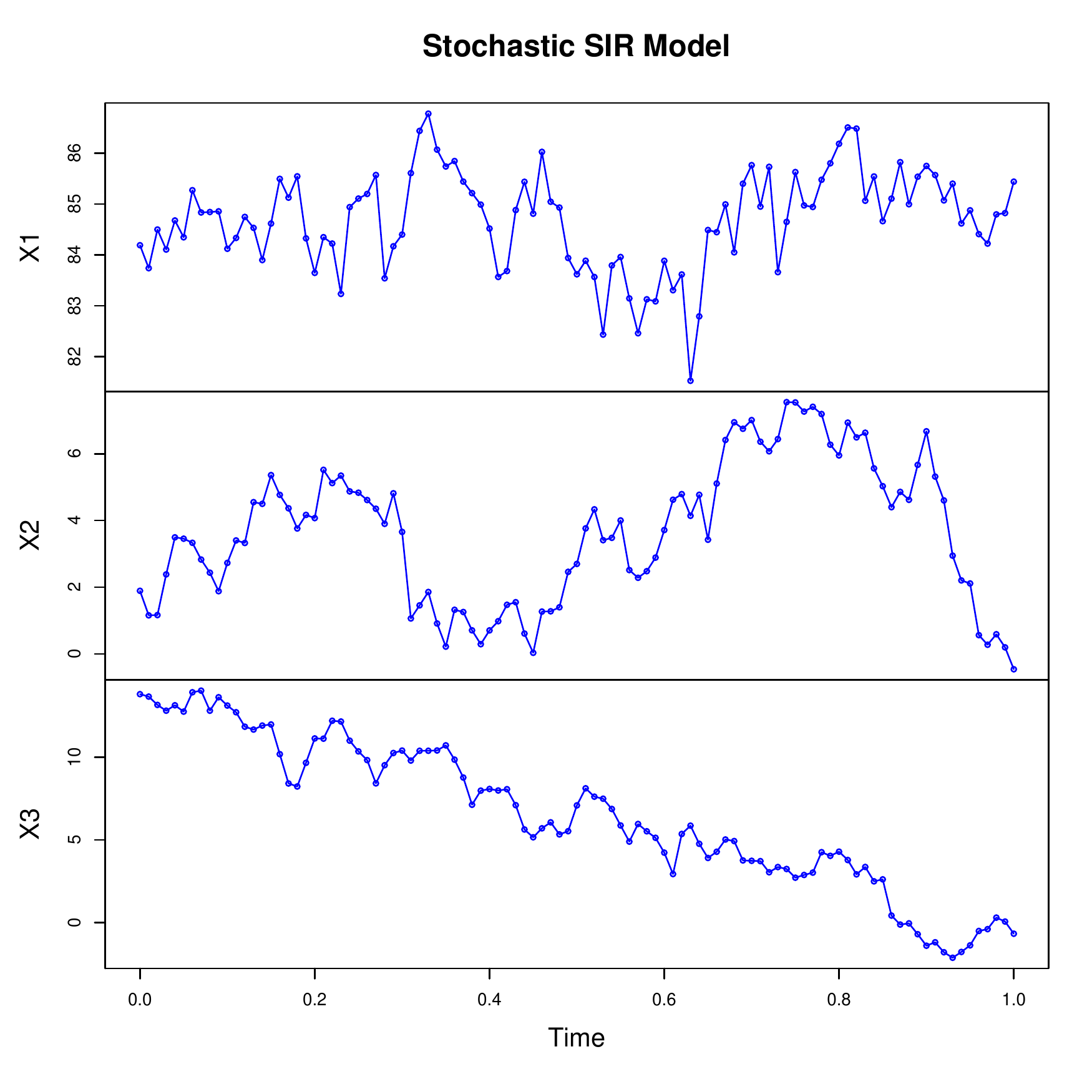}
		\caption{SIR Model with magnified S($X1$), I($X2$) and R($X3$) with $\sigma_1=0.05$, $\sigma_2=0.08557$ and $\sigma_3=0.12$.}
		\label{fig11}
	\end{minipage}
\end{figure}

\begin{figure}[H]
	\begin{minipage}[b]{0.4\linewidth}
		\centering
		\includegraphics[width=\textwidth]{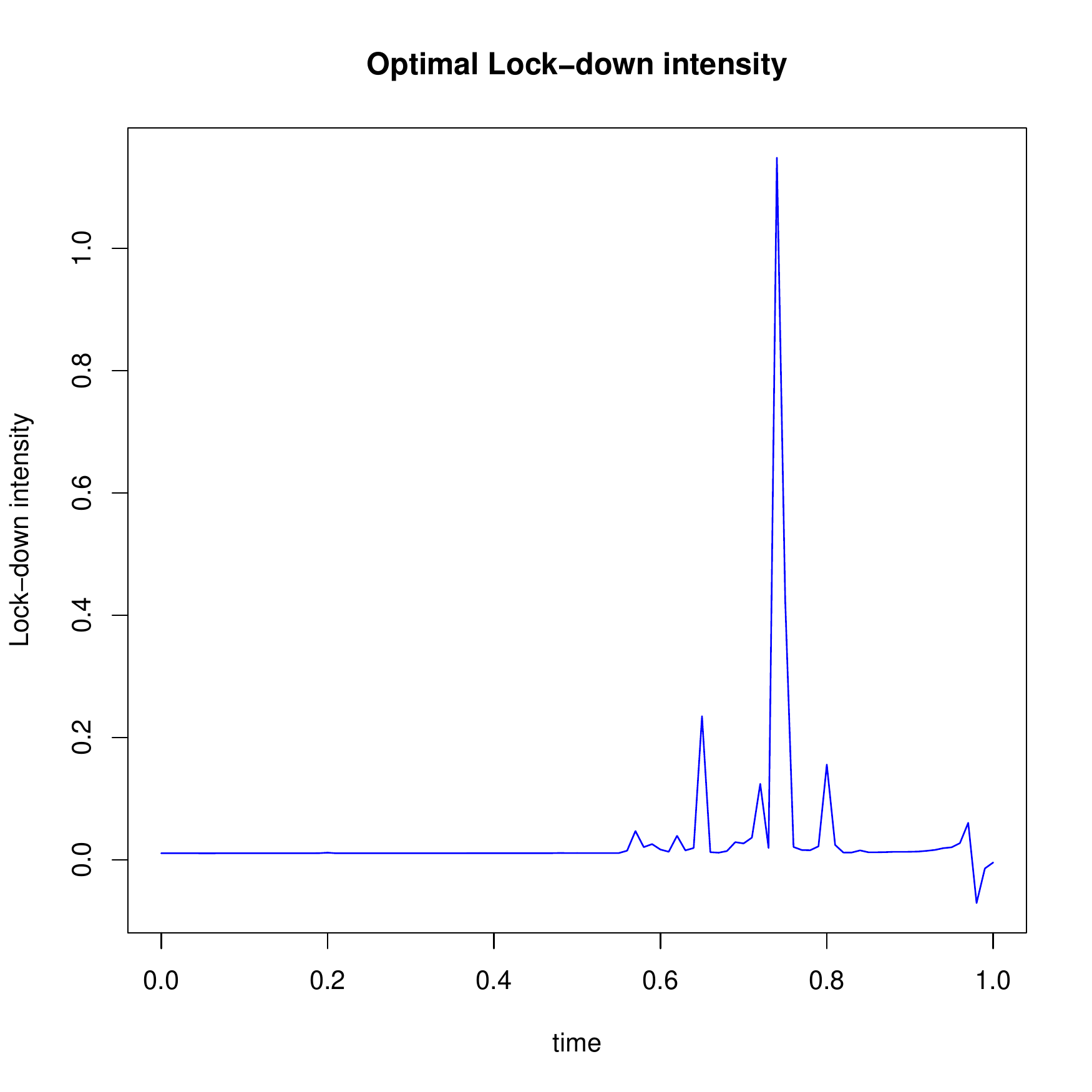}
		\caption{Optimal lock-down of UK data for first $100$ days of $2021$ with diffusion coefficient $\sigma_2=0.08557$.}
		\label{fig12}
	\end{minipage}
	\hspace{0.5cm}
	\begin{minipage}[b]{0.4\linewidth}
		\centering
		\includegraphics[width=\textwidth]{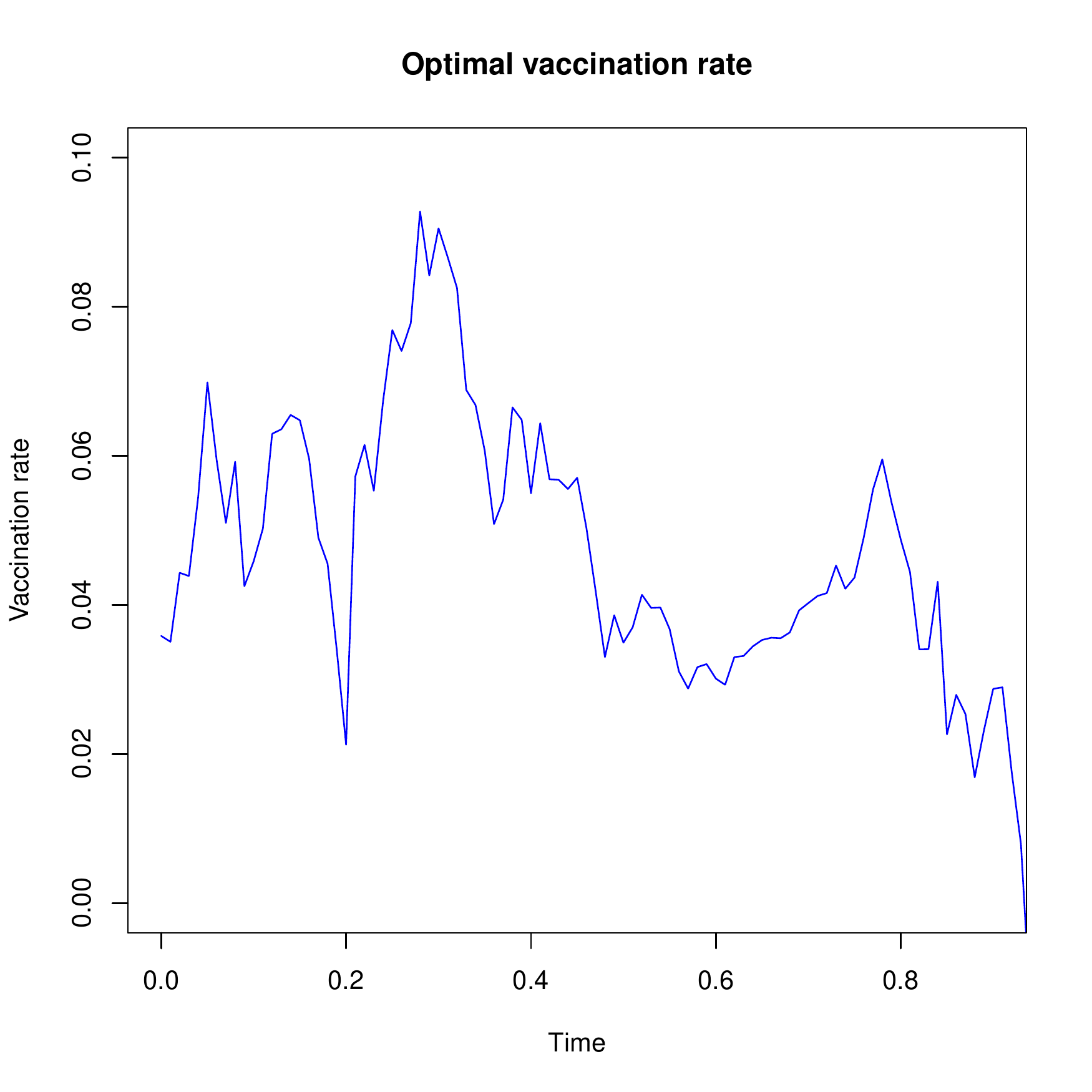}
		\caption{Optimal vaccination rate of UK data for first $100$ days of $2021$ with diffusion coefficient $\sigma_2=0.08557$.}
		\label{fig13}
	\end{minipage}
\end{figure}

To generate figures \ref{fig10}-\ref{fig13} we consider first $100$ days of $2021$ starting from January 1. Then we convert this time interval between $0$ and $1$. Figures \ref{fig10}-\ref{fig11} describe the stochastic SIR model at the beginning of 2021. Furthermore, figure \ref{fig11} magnifies each of stochastic susceptible ($X1$), infected ($X2$) and recovered ($X3$) curves. $X2$ curve does not have any pattern because of the high volatility of UK data (i.e. $\sigma_2=0.08557$) which is consistent with the simulation in figure \ref{fig7}. Although the susceptible curve has a downward trend in both figures \ref{fig6} and \ref{fig10}, further magnification of the behavior leads us to more volatility as explained by curve $X1$ in figure \ref{fig11}. The stochastic recovery curve shows a similar pattern to the theoretical results. 

\medskip

Optimal lock-down intensity curves in figures \ref{fig8} and \ref{fig12} show similar trend between days $1$ and $50$. After $50^{th}$ day the intensity curve in figure \ref{fig12} becomes more ergodic and shows a spike between days $70$ and $80$. This tells us that although every agent in the COVID-19 environment has perfect and complete information about the pandemic, after a certain point in time they want to leave the house probably to purchase necessary items or just because of some social interactions. Contrarily, optimal vaccination rate curves presented in figures \ref{fig9} and \ref{fig13} are quite different in character. The curve in figure \ref{fig9} shows an upward trend in vaccination rate, while the curve in figure \ref{fig13} does not show any trend. As we take the UK data at the beginning of 2021, the infection and death rates have a huge spike. Furthermore, since at that time new vaccines are coming slowly into the economy, people have less confidence in them, and this leads to an unstable trend toward optimal vaccination rate.

\medskip

\section{Conclusion}

\medskip

In this paper, a stochastic pandemic SIR model with a non-linear incidence rate 
$\be SI/(1+\rho I+\eta N)$ and a stochastic dynamic infection rate is considered. We use a Feynman-type path integral approach to obtain optimal lock-down intensity and vaccination rate because simulating an HJB equation is almost impossible due to the curse of dimensionality \citep{polanskyrus,pr5,pr6}. The main aspect of this study lies in the aspect of the existence of global stability and uniqueness of the control variables when the information is \emph{perfect} and \emph{complete}. Furthermore, we can show the existence of a contraction mapping point in the Brouwer sense.

To determine the infection dynamics, we have divided the immunity level into five subcategories such as \emph{very low, somewhat low, medium, somewhat high} and \emph{very high}. We have used Erdos-Renyi random graph model to investigate the infection rate among agents with different levels. We have minimized an agent's cost of COVID-19 subject to a stochastic SIR and infection dynamics. By utilizing a Feynman-type path integral approach we determine a Fokker-Plank type equation and obtain an optimal lock-down intensity and vaccination rate. We also did some simulation studies based on the parameters in \cite{caulkins2021optimal} and \cite{rao2014}. Since we assume all agents in the pandemic environment are risk averse, the optimal lock-down intensity went up at the beginning of our time interval and then came very close to zero \citep{ahamed2021i,ahamed2021,ahamed2021det,ahamed2022}. The reason behind this is that due to the availability of \emph{perfect} and \emph{complete} information an individual does not want to go out and gets infected by this pandemic. At the end of our time of the study, we observe that lock-down intensity is slightly improved although the optimal vaccination rate has increased over the time interval we have studied.

Data analysis tells us that, although people are risk averse to COVID-19, after a certain point of time they come out of their homes to do social interaction probably because of the necessity to purchase food or some other important items. As of the beginning of 2021, the incidence of COVID-19 has experienced a spike, with the new incidence of vaccines people have less faith in medicines. This leads to an unprecedented movement of \emph{optimal vaccination rate} in the first $100$ days of that year. 

\medskip

\section*{Appendix}

\subsection*{Proof of Proposition \ref{p0}}

\medskip

For $M[\mathbf X(s)]\in(0,\infty)$ let $\{\mathbf X_n\}_{n\geq 1}$ be a minimizing sequence in the convex set $\mathcal X$ so that,
\begin{equation*}
\lim_{n\ra\infty}\E_0\left\{\int_0^tc[u(s),\mathbf X_n(s)]ds\bigg|\mathcal F_0\right\}=M[\mathbf X(s)]>0.
\end{equation*}
By compactness theorem A.3.5 in \cite{pham2009} there exists a convex combination ${\widetilde{\mathbf X}_n(s)}\in \text{conv}\left(\mathbf X_n,\mathbf X_{n+1},...\right)\in\mathcal X$ so that ${\widetilde{\mathbf X}_n(s)}\overset{a.s.}\to \mathbf X^*$. As $\{\mathbf X_n\}_{n\geq 1}$ is a minimizing sequence hence, $\mathbf X^*\in\mathcal X$. Convexity of the cost function implies
\[
\lim_{n\ra\infty}\E_0\left\{\int_0^tc[u(s),{\widetilde{\mathbf X}_n(s)}]ds\bigg|\mathcal F_0\right\}=M[\mathbf X(s)]>0.
\]
Above condition and convexity of $c$ implies $\inf_n\E_0\left\{\int_0^tc[u(s),{\widetilde{\mathbf X}_n(s)}]ds\bigg|\mathcal F_0\right\}>0$. Therefore, optimality of $\mathbf X(s)$ i.e.
\[
\E_0\left\{\int_0^tc[u(s),\mathbf X^*(s)]ds\bigg|\mathcal F_0\right\}=M[\mathbf X(s)],
\]
is achieved iff we are able to show
\begin{equation}\label{2}
\lim_{n\ra\infty}\E_0\left\{\int_0^tc[u(s),{\widetilde{\mathbf X}_n(s)}]ds\bigg|\mathcal F_0\right\}=\E_0\left\{\int_0^tc[u(s),\mathbf X^*(s)]ds\bigg|\mathcal F_0\right\},
\end{equation}
which is uniform integrability of $\left\{c[u(s),\mathbf X_n(s)]\right\}_{n\geq 1}$. If $c[u(s),\infty]>0$ and define the initial condition
\[
\mathbf X_0:=\inf\left\{\mathbf X(s)>0:c[u(s),\mathbf X(s)]>0\right\}<\infty.
\]
We will prove by contradiction by assuming that sequence $\left\{c[u(s),\mathbf X_n(s)]\right\}_{n\geq 1}$ is not uniformly integrable. Hence, $\exists\ \delta>0$ so that.
\[
\lim_{n\ra\infty}\E_0\left\{\int_0^tc[u(s),{\widetilde{\mathbf X}_n(s)}]ds\bigg|\mathcal F_0\right\}=\E_0\left\{\int_0^tc[u(s),\mathbf X^*(s)]ds\bigg|\mathcal F_0\right\}+2\delta.
\]
A subsequence $\{{\widetilde{\mathbf X}_n(s)}\}_{n\geq 1}$ and Corollary A.1.1 in \cite{pham2009} implies that there exist disjoint sets $(E_n)_{n\geq 1}$ in $(\mathcal X_\infty,\mathcal F_0,\mathcal P)$ so that,
\[
\E_0\left\{\int_0^tc[u(s),{\widetilde{\mathbf X}_n(s)}]\mathbbm{1}_{E_n}ds\bigg|\mathcal F_0\right\}\leq\delta,\ \text{for all}\ n\geq 1.
\]
There exists a sequence of state variables in $(\mathcal X_\infty,\mathcal F_0,\mathcal P)$
\[
J_n(s)=\mathbf X_0+\sum_{l=1}^n{\widetilde{\mathbf X}_l(s)}\mathbbm{1}_{E_n}.
\]
For any local probability measure $\mathcal Q$ with local martingale $\mathcal M(\mathcal X)$ with the relationship $\mathcal Q\in\mathcal M(\mathcal X)$ yields,
\[
\E_0^{\mathcal Q}\left[J_n(s)\bigg|\mathcal F_0\right]\leq\mathbf X_0+\sum_{l=1}^n\E_0^{\mathcal Q}\left[{\widetilde{\mathbf X}_l(s)}\bigg|\mathcal F_0\right]\mathbbm{1}_{E_n}\leq\mathbf X_0+n\mathbf X,
\]
as ${\widetilde{\mathbf X}_l(s)}\in\mathcal X$. Clearly, $J_n(s)\in\mathcal X(\mathbf X_0+n\mathbf X)$ where $\mathcal X(\mathbf X_0+n\mathbf X)$ stands for a convex functional space of state variables such that the property $\mathbf X_0+n\mathbf X$ holds. Therefore,
\begin{multline*}
\E_0\left\{\int_0^tc[u(s),J_n(s)]ds\bigg|\mathcal F_0\right\}=\E_0\left\{\int_0^tc\left[u(s),\mathbf X_0+\sum_{l=1}^n{\widetilde{\mathbf X}_l(s)}\mathbbm{1}_{E_n}\right]ds\bigg|\mathcal F_0\right\}\\
\leq\E_0\left\{\int_0^tc\left[u(s),\sum_{l=1}^n{\widetilde{\mathbf X}_l(s)}\mathbbm{1}_{E_n}\right]ds\bigg|\mathcal F_0\right\}=\sum_{l=1}^n\E_0\left\{\int_0^tc\left[u(s),{\widetilde{\mathbf X}_l(s)}\mathbbm{1}_{E_n}\right]ds\bigg|\mathcal F_0\right\}\leq \delta n.
\end{multline*}
By convexity of the dynamic cost function we get,
\begin{equation*}
\liminf_{\mathbf X\ra\infty}\frac{M[\mathbf X(s)]}{\mathbf X(s)}\leq \liminf_{\mathbf X\ra\infty}\frac{\E_0\left\{\int_0^tc[u(s),J_n(s)]ds\bigg|\mathcal F_0\right\}}{\mathbf X_0+n\mathbf X}\leq\liminf_{\mathbf X\ra\infty}\frac{n\delta}{\mathbf X_0+n\mathbf X}\leq\delta.
\end{equation*}
After carefully setting $\delta\ra 0$ we conclude that $\liminf_{\mathbf X\ra\infty}M[\mathbf X(s)]/\mathbf X(s)< 0$ which is a contradiction from the assumption $\liminf_{\mathbf X\ra\infty}M[\mathbf X(s)]/\mathbf X(s)\geq 0$. Therefore, condition explained in Equation (\ref{2})is true and $\mathbf X^*$ is the solution to $M[\mathbf X(s)]$ for all $s\in[0,t]$. The uniqueness follows from the strict convexity of the cost function on $(0,\infty)$ and known filtration process $\mathcal F_0$. $\square$

\subsection*{Proof of Proposition \ref{p1}}

\medskip

For small continuous time interval $[s,\tau]$ It\^o formula yields 
\begin{multline*}
\mathcal Z(\tau,\mathbf X(\tau))=\mathcal Z(s,\mathbf X(s))\\
+\int_s^\tau\left\{\frac{\partial}{\partial s}\mathcal Z(\nu,\mathbf X)	+\bm\mu(\nu,u,\mathbf X)\frac{\partial}{\partial\mathbf X}\mathcal Z(\nu,\mathbf X)+\frac{1}{2}\text{$trace$}\left\{\bm\sigma^T(\nu,\mathbf X)\left[\frac{\partial^2}{\partial X^T \partial X} \mathcal Z(\nu,\mathbf X)\right]\bm\sigma(\nu,\mathbf X)\right\}\right\}d\nu\\
+\int_s^\tau \bm\sigma(s,\mathbf X)\frac{\partial}{\partial \mathbf X}\mathcal Z(s,\mathbf X) d\mathbf B(\nu).
\end{multline*}
As we have already assumed that $\bm\sigma(s,\mathbf X)\frac{\partial}{\partial \mathbf X}\mathcal Z(s,\mathbf X)$ is in Hilbert space $L^2$, then in the above equation integral part with respect to time vanishes \citep{lindstrom2018}. Applying boundary condition $\mathcal Z(\tau,\mathbf X)=\Phi(\mathbf X)$, the initial condition $\mathbf X(s)=\mathbf X_s$ and taking conditional expectation on the remainder part of the above equation yield 
\[
\mathbf Z(s,\mathbf X_s)=\E_s[\Phi(\mathbf X(\tau))].
\] 
This completes the proof.     $\square$

\medskip

\subsection*{Proof of Proposition \ref{p2}}

\medskip

Suppose $\mathbf X$ and $\widetilde{\mathbf X}$ both are strong solutions on the $4$-dimensional Brownian motion $\mathbf B(s)$ for all $s\in[0,t]$ under complete probability space $\{\mathcal X_\infty,\mathcal F_0,\mathcal P\}$. Define stopping times 
\begin{align*}
s_\rho&:=\inf\left\{s\geq 0;||\mathbf X(s)||\geq\rho,\ \forall\ \rho\geq 1\right\}\\
\tilde{s}_\rho&:=\inf\left\{s\geq 0;||\widetilde{\mathbf X}(s)||\geq\rho,\ \forall\ \rho\geq 1\right\}.
\end{align*}
Setting $\mathcal S_\rho\triangleq s_\rho\wedge\tilde{s}_\rho$ yields $P\left[\lim_{\rho\ra\infty}\mathcal S_\rho\right]\overset{a.s.}\to\infty$ and
\begin{multline*}
\mathbf X(s\wedge\mathcal S_\rho)-\widetilde{\mathbf X}(s\wedge\mathcal S_\rho)=\int_0^{s\wedge\mathcal S_\rho}\left\{\bm\mu\left[\nu,u(\nu),\mathbf X(\nu)\right]-\bm\mu\left[\nu,u(\nu),\widetilde{\mathbf X}(\nu)\right]\right\}d\nu\\
+\int_0^{s\wedge\mathcal S_\rho}\left\{\bm\sigma\left[\nu,\mathbf X(\nu)\right]-\bm\sigma\left[\nu,\tilde{\mathbf X}(\nu)\right]\right\}d\mathbf B(\nu)
\end{multline*}

For any finite constant $\mathcal K$, H\"older inequality for Lebesgue integrals, property $3.2.27$ of \cite{karatzas2012} and Assumption \ref{a1} imply
\begin{multline*}
\E_0\left\{\left|\left|\mathbf X(s\wedge\mathcal S_\rho)-\widetilde{\mathbf X}(s\wedge\mathcal S_\rho)\right|\right|^2\biggr|\mathcal F_0\right\}\\
\leq 9\E_0\left\{\left[\int_0^{s\wedge\mathcal S_\rho}\left|\left|\bm\mu\left[\nu,u(\nu),\mathbf X(\nu)\right]-\bm\mu\left[\nu,u(\nu),\widetilde{\mathbf X}(\nu)\right]\right|\right|d\nu\right]^2\biggr|\mathcal F_0\right\}\\
+9\E_0\left\{\sum_{k=1}^d\left[\sum_{l=1}^r\int_0^{s\wedge\mathcal S_\rho}\left[\bm\sigma_{kl}\left[\nu,\mathbf X(\nu)\right]-\bm\sigma_{kl}\left[\nu,\tilde{\mathbf X}(\nu)\right]\right]d\mathbf B^{(l)}(\nu)\right]^2\biggr|\mathcal F_0\right\}\\
\leq 9s\E_0\left\{\int_0^{s\wedge\mathcal S_\rho}\left|\left|\bm\mu\left[\nu,u(\nu),\mathbf X(\nu)\right]-\bm\mu\left[\nu,u(\nu),\widetilde{\mathbf X}(\nu)\right]\right|\right|^2d\nu\biggr|\mathcal F_0\right\}\\
+9\E_0\left\{\int_0^{s\wedge\mathcal S_\rho}\left|\left|\bm\sigma\left[\nu,\mathbf X(\nu)\right]-\bm\sigma\left[\nu,\widetilde{\mathbf X}(\nu)\right]\right|\right|^2d\nu\biggr|\mathcal F_0\right\}\\
\leq 9(1+t)\mathcal K^2\int_0^t\E_0\left\{\left|\left|\mathbf X(s\wedge\mathcal S_\rho)-\widetilde{\mathbf X}(s\wedge\mathcal S_\rho)\right|\right|^2d\nu\biggr|\mathcal F_0\right\}.
\end{multline*}
Following \cite{karatzas2012} we know for $s\in[0,t]$ above condition implies
\begin{equation*}
\E_0\left\{\left|\left|\mathbf X(s\wedge\mathcal S_\rho)-\widetilde{\mathbf X}(s\wedge\mathcal S_\rho)\right|\right|^2\biggr|\mathcal F_0\right\}\leq 9t\mathcal K^2+9\mathcal K^2\int_0^t\E_0\left\{\left|\left|\mathbf X(s\wedge\mathcal S_\rho)-\widetilde{\mathbf X}(s\wedge\mathcal S_\rho)\right|\right|^2d\nu\biggr|\mathcal F_0\right\}.
\end{equation*}
Therefore, $\{\mathbf X(s\wedge\mathcal S_\rho);s\in[0,t]\}$ and $\widetilde{\mathbf X}(s\wedge\mathcal S_\rho);s\in[0,t]$ are modification of each other and hence, are indistinguishable. Allowing $\rho\ra\infty$ gives us $\{\mathbf X(s);s\in[0,t]\}$ and $\{\widetilde{\mathbf X}(s);s\in[0,t]\}$ are indistinguishable.  $\square$

\subsection*{Proof of Proposition \ref{p3}}

\medskip

For each optimal solution ${\mathbf X}^*\in\mathbb F^2$ of Equation (\ref{5}), define a squared integrable progressively measurable process $Z(\mathbf X^*)$ by 
\begin{equation}\label{7}
Z(\mathbf X^*)_s=\mathbf X(0)+\int_0^t\bm\mu(s,u,\mathbf X)ds+\int_0^t\bm\sigma(s,\mathbf X)d\mathbf B(s).
\end{equation} 
We will show that $Z(\mathbf X^*)\in\mathbb F^2$. Furthermore, as $\mathbf X^*$ is a solution of Equation (\ref{5}) iff $Z(\mathbf X^*)=\mathbf X^*$, we will show that $Z$ is the strict contraction of the Hilbert space $\mathbb F^2$. Using the fact that 
\[
|\bm\mu(s,u,\mathbf X)|^2\leq c_0\left[1+|\mathbf X|^2+|\bm\mu(s,u,\mathbf X(0))|^2|\right]
\]
yields
\begin{equation}\label{8}
||Z(\mathbf X)||^2\leq 4\left[t \E|\mathbf X(0)|^2+\E\int_0^t\bigg|\int_0^s\bm\mu(s',u,\mathbf X)ds'\bigg|^2ds+t\E\sup_{0\leq s\leq t}\bigg|\int_0^s\bm\sigma(s',\mathbf X(s'))d\mathbf B(s')\bigg|^2ds\right].
\end{equation}
Assumption \ref{a2} implies $t\E|\mathbf X(0)|^2<\infty$. It will be shown that the second and third terms of the right hand side of the inequality (\ref{6}) are also finite. Assumption \ref{a1} implies,
\begin{multline*}
\E\int_0^t\bigg|\int_0^s\bm\mu(s',u,\mathbf X)ds'\bigg|^2ds\leq \E\int_0^ts\left(\int_0^s|\bm\mu(s',u,\mathbf X)|^2ds'\right)ds\\\leq c_0\E\int_0^t s\left(\int_0^s(1+|\bm\mu(s',u,\mathbf X(0))|^2+|\mathbf X(s)|^2)ds'\right)ds\\\leq c_0t^2\left(1+||\bm\mu(s',u,\mathbf X(0))||^2+\E\sup_{0\leq s\leq t}|\mathbf X(s)|^2\right)<\infty.
\end{multline*}
Doob's maximal inequality and Lipschitz assumption (i.e. Assumption \ref{a1}) implies,
\begin{multline*}
t\E\sup_{0\leq s\leq t}\bigg|\int_0^s\bm\sigma(s',\mathbf X(s'))d\mathbf B(s')\bigg|^2ds\leq 4t\E\int_0^t|\bm\sigma(s,\mathbf X(s))|^2ds\\\leq 4c_0\E\int_0^t(1+|\bm\sigma(\mathbf X(0))|^2+|\mathbf X(s)|^2)ds\\\leq 4c_0t^2\left(1+||\bm\sigma(\mathbf X(0))||^2+\E\sup_{0\leq s\leq t}|\mathbf X(s)|^2\right)<\infty.
\end{multline*}
As $Z$ maps $\mathbb F^2$ into itself, we show that it is strict contraction. To do so we change Hilbert norm $\mathbb F^2$ to an equivalent norm. Following \cite{carmona2016} for $a>0$ define a norm on $\mathbb F^2$ by $$||\xi||_a^2=\E\int_0^t\exp(-as)|\xi_s|^sds.$$ If $\mathbf X(s)$ and $\mathbf Y(s)$ are generic elements of $\mathbb{F}^2$ where $\mathbf X(0)=\mathbf Y(0)$, then
\begin{multline*}
\E|Z(\mathbf X(s))-z(\mathbf Y(s))|^2\leq 2\E\bigg|\int_0^\tau[\bm\mu(s',u,\mathbf X(s'))-\bm\mu(s',u,\mathbf Y(s'))]ds\bigg|^2\\+2\E\bigg|\int_0^\tau[\sigma_0^k(\mathbf X(s'))-\bm\sigma(s',\mathbf Y(s'))]d\mathbf B(s')\bigg|^2\\\leq 2\tau\E\int_0^\tau|\bm\mu(s',u,\mathbf X(s'))-\bm\mu(s',u,\mathbf Y(s'))|^2ds'+2\E\int_0^\tau|\bm\sigma(s',\mathbf X(s'))-\bm\sigma(s',\mathbf Y(s'))|^2ds'\\\leq c_0(1+\tau)\int_0^\tau\E|\mathbf X(s')-\mathbf Y(s')|^2ds',
\end{multline*}
by Lipschitz's properties of drift and diffusion coefficients. Hence.
\begin{multline*}
||Z(\mathbf X)-Z(\mathbf Y)||_a^2=\int_0^t\exp(-as)\E|Z(\mathbf X(s)-Z(\mathbf Y(s)))|^2ds\\
\leq c_0t\int_0^t\exp(-as)\int_0^t\E|\mathbf X(s')-\mathbf Y(s')|^2ds'ds\\\leq c_0t\int_0^t\exp(-as)ds\int_0^t\E|\mathbf X(s')-\mathbf Y(s')|^2ds'\leq\frac{c_0t}{a}||\mathbf X-\mathbf Y||_a^2.
\end{multline*}
Furthermore, if $c_0t$ is very large, $Z$ becomes a strict contraction. Finally, for $s\in[0,t]$
\begin{multline*}
\E \sup_{0\leq s\leq t}|\mathbf X(s)|^2=\E \sup_{0\leq s\leq t}\bigg|\mathbf X(0)+\int_0^{s'}\bm\mu(r,u,\mathbf X(r))dr+\int_0^{s'}\bm\sigma(r,\mathbf X(r))d\mathbf B(r)\bigg|^2\\\leq 4\left[\E|\mathbf X(0)|^2+s\E\int_0^s|\bm\mu(s',u,\mathbf X(s'))|^2ds'+4\E\int_0^s|\bm\sigma(s',\mathbf X(s'))|ds'\right]\\\leq c_0\left[1+\E|\mathbf X(0)|^2+\int_0^s\E\sup_{0\leq r\leq s'}|\mathbf X(r)|^2dr\right],
\end{multline*}
where the constant $c_0$ depends on $t$, $||\bm\mu||^2$ and $||\bm\sigma||^2$. Gronwall's inequality implies,
\begin{equation*}
\E \sup_{0\leq s\leq t}|\mathbf X(s)|^2\leq c_0(1+\E|\mathbf X(0)|^2)\exp{(c_0t)}.\ \ \square
\end{equation*}

\medskip

\subsection*{Proof of Proposition \ref{p4}}

\medskip

In order to show Condition \ref{p4.0} we will use Hahn-Jordan orthogonal decomposition of total variation \citep{del2004}
\[
\mathcal P=\mathcal P_1-\mathcal P_2=\mathcal P^+-\mathcal P^-,
\]
such that $||\mathcal P_1-\mathcal P_2||_{tv}=\mathcal P^+(H)=\mathcal P^-(H)$. Therefore, for quantum Lagrangian $\mathcal L\in H$ we have
\begin{multline*}
|\mathcal P_1(\mathcal L)-\mathcal P_2(\mathcal L)|=\left|\int_{\mathbf R^4}\mathcal L(s,u,\mathbf X)\mathcal P^+(d\mathbf X)-\int_{\mathbf R^4}\mathcal L(s,u,\widetilde{\mathbf X})\mathcal P^-(d\widetilde{\mathbf X})\right|\\
=||\mathcal P_1-\mathcal P_2||_{tv}\left|\int_{\mathbf R^4}\left[\mathcal L(s,u,\mathbf X)-\mathcal L(s,u,\widetilde{\mathbf X})\right]\frac{\mathcal P^+(d\mathbf X)}{\mathcal P^+(H)}\frac{\mathcal P^-(d\widetilde{\mathbf X})}{\mathcal P^-(H)}\right|.
\end{multline*}
Above condition implies,
\[
|\mathcal P_1(\mathcal L)-\mathcal P_2(\mathcal L)|\leq ||\mathcal P_1-\mathcal P_2||_{tv}.
\]
Supremum over all $\mathcal L\in H$ yields,
\[
\sup\left\{|\mathcal P_1(\mathcal L)-\mathcal P_2(\mathcal L)|;\ \mathcal L\in H\right\}\leq ||\mathcal P_1-\mathcal P_2||_{tv}.
\]
The reverse inequality can be checked by using the simple function $\mathbbm 1_B$ such as $B\in\mathcal E\in H$. Now we will show Condition \ref{p4.1}. By the construction of this pandemic framework there exist two non-interacting neighborhoods of agents $H_+$ and $H_-$ so that $\mathcal P^+(H_-)=0=\mathcal P^-(H_+)$. Hence, for all $B\in\mathcal E$, we have  
\begin{equation*}
\mathcal P^+(B)=\mathcal P(B\cap H_+)\geq 0\ \ \text{and}\ \ \mathcal P^-(B)=-\mathcal P^-(B\cap H_-)\geq 0,
\end{equation*} 
which implies, 
\begin{equation}\label{10}
\mathcal P_1(B\cap H_+)\geq\mathcal P_2(B\cap H_+)\ \ \text{and}\ \ \mathcal P_2(B\cap H_-)\geq\mathcal P_1(B\cap H_-).
\end{equation}
For any $B\in\mathcal E$ define $\eta$ as $\eta\triangleq \mathcal P_1(B\cap H_-)+\mathcal P_2(B\cap H_+)$. By construction of this pandemic network we have
\begin{equation}\label{11}
\eta(B)\leq \mathcal P_1(B)\wedge\mathcal P_2(B)\ \ \text{and}\ \ \eta(H)=\mathcal P_1(H_-)+\mathcal P_2(H_+).
\end{equation}
As total variation distance between two immunity group is
\[
||\mathcal P_1-\mathcal P_2||_{tv}=\mathcal P^+(H)=\mathcal P(H_+)=\mathcal P_1(H_+)-\mathcal P_2(H_-)=1-\left[\mathcal P_1(H_+)+\mathcal P_2(H_-)\right],
\]
Condition \ref{11} implies
\[
1-\sup_{\gamma\leq\mathcal P_1,\mathcal P_2}\gamma(H)\leq1-\eta(H)=||\mathcal P_1-\mathcal P_2||_{tv}.
\] 
In order to show the reverse inequality assume $\gamma$ be a non-negative such as for all $B\in\mathcal E$ yields $\gamma(B)\leq\mathcal P_1(B)\wedge\mathcal P_2(B)$. Suppose, if we consider $B=H_+$ and $B=H_-$ respectively, then we have
\begin{equation*}
\gamma(H_+)\leq\mathcal P_1(H_+)\ \ \text{and}\ \ \gamma(H_-)\leq\mathcal P_2(H_-),
\end{equation*}
which yields
\begin{equation*}
\gamma(H)\leq \mathcal P_1(H_+)+\mathcal P_2(H_-)=1-||\mathcal P_1-\mathcal P_2||_{tv}.
\end{equation*}
Therefore, $1-\gamma(H)\geq ||\mathcal P_1-\mathcal P_2||_{tv}$. Finally, taking the infimum of over all the distributions of $\gamma\leq\mathcal P_1$ and $\mathcal P_2$, Condition \ref{p4.1} is obtained. To show Condition \ref{p4.2} we are going to use the similar idea like above. First, by using \ref{10} define $\mathcal P_2(H_+)=\mathcal P_1(H_+)\wedge\mathcal P_2(H_+)$ and $\mathcal P_1(H_-)=\mathcal P_1(H_-)\wedge\mathcal P_2(H_-)$. This yields
\begin{equation*}
\eta(H)=\mathcal P_1(H_-)+\mathcal P_2(H_+)=\left[\mathcal P_1(H_-)\wedge\mathcal P_2(H_-)\right]+\left[\mathcal P_1(H_+)\wedge\mathcal P_2(H_+)\right].
\end{equation*}
Non-interaction of agents between $H_+$ and $H_-$ implies
\[
\eta(H)\leq\inf\sum_{k=1}^K\left[\mathcal P_1(B_k)\wedge\mathcal P_2(B_k)\right],
\]
where the infimum is taken over all resolutions of $H$ into pairs of non-interacting subgroups $B_k, \ k\in[1,K],\ K\geq 1$. To show the reverse inequality we use the definition of $\eta$ \citep{del2004}. Using Condition \ref{11} we know for any finite resolution $B_k\in\mathcal E$ the inequality $\eta(B_k)\leq\mathcal P_1(B_k)\wedge\mathcal P_2(B_k)$ holds. Thus,
\[
\eta(H)=\sum_{k=1}^K\eta(B_k)\leq \sum_{k=1}^K\mathcal P_1(B_k)\wedge\mathcal P_2(B_k).
\] 
Taking the infimum of all resolutions and using $\eta(H)=1-||\mathcal P_1-\mathcal P_2||_{tv}$ yield Condition \ref{p4.2}. This completes the proof. $\square$

\medskip

\subsection*{Proof of Proposition \ref{p5}}

\medskip

We have divided the proof into two cases.

$\mathbf{Case\ I}$: We assume that $m\subset \mathbb N$, a set $\beth$ with condition $|\beth|=m+1$, and affinely independent state variables, vaccination rates and lock-down intensities $\{\mathbf X_k(s)\}_{k\in \beth}\subset \mathbb R^{6}$ such that $\widetilde \Xi$ coincides with the simplex convex set of $\{ \mathbf X_k(s)\}_{k\in \beth}$. For each $Z(s)\subset \Xi$, there is a unique way in which the vector $Z(s)$ can be written as a convex combination of the extreme valued state variables, vaccination rates and lock-down intensities; such as, $Z(s)=\sum_{k\in\beth}\alpha_k(s,\mathbf X)\mathbf X_k(s)$ so that $\sum_{k\in\beth}\alpha_k(s,\mathbf X)=1$ and $\alpha_k(s,\mathbf X)\geq 0,\ \forall k\in\beth$ and $s\in[0,t]$. For each $k\in\beth$, define a set 
\[
\widetilde\Xi_k\triangleq\left\{Z\in\widetilde\Xi:\a_k[\mathcal L(s,u,\mathbf X)]\leq\a_k(s,\mathbf X)\right\}.
\]
By the continuity of the quantum Lagrangian of a agent $\mathcal L$, for each $k\in\beth$, $\widetilde\Xi_k$ is closed. Now we claim that, for every $\widetilde \beth \subset \beth$, the convex set consists of $\{\mathbf X_k\}_{k\in\widetilde\beth}$ is proper subset of $\bigcup_{k\in\tilde\beth}\widetilde\Xi_k$. Suppose $\widetilde\beth\subset\beth$ and $Z(s)$ is also in the non-empty, convex set consists of the state variables, vaccination rates and the lock-down intensities $\{\mathbf X_k(s)\}_{k\in\tilde\beth}$. Thus, $\sum_{k\in\widetilde\beth}\a_k(s,\mathbf X)=1\geq\sum_{k\in\widetilde\beth}\a_k[\mathcal L(s,u,\mathbf X)]$. Therefore, there exists $k\in\widetilde\beth$ such that $\a_k(s,\mathbf X)\geq\a_k\left[\mathcal L(s,u,\mathbf X)\right]$ which implies $Z(s)\in\tilde\Xi\subset \bigcup_{l\in\tilde\beth}\tilde\Xi_l$. By \emph{Knaster-Kuratowski-Mazurkiewicz Theorem}, there is $\bar {\mathbf X}_k^*\in \bigcap_{k\in\beth}\widetilde\Xi_k$, in other words, the condition $\a_k\left[\mathcal L(s,u^*,\bar {\mathbf X}_k^*)\right]\leq\a_k(s,\bar {\mathbf X}_k^*)$ for all $k\in\beth$ and for each $s\in[0,t]$ \citep{gonzalez2010}. Hence, $\mathcal L(s,u^*,\bar {\mathbf X}_k^*)=\bar {\mathbf X}_k^*$ or $\mathcal L$ has a fixed-point. 

$\mathbf{Case\ II}$: Again consider $\widetilde\Xi\subset \mathbb R^6$ is a non-empty, convex and compact set. Then for $m\subset \mathbb N$, a set $\beth$ with condition $|\beth|=m+1$, and affinely independent state variables, vaccination rates and lock-down intensities $\{\mathbf X_k(s)\}_{k\in \beth}\subset \mathbb R^{6}$ such that $\widetilde\Xi$ is a proper subset of the convex set based on $\{\mathbf X_k(s)\}_{k\in\beth}$ for all $s\in[0,t]$. Among all the simplices, suppose $\hat\aleph$ is the set with smallest $m$. Let $\tilde Z(s)$ be a dynamic point in the $m$-dimensional interior of $\hat\aleph$. Define ${\hat{\mathcal L}}$, an extension of $\mathcal L$ to the whole simplex $\hat\aleph$, as follows. For every $Z(s)\in\hat\aleph$, let 
\[
\bar\zeta(s,Z):\max\left\{\bar\zeta\in[0,1]:(1-\bar\zeta)\tilde Z(s)+\bar\zeta Z(s)\in\widetilde\Xi\right\},\ \forall s\in[0,1],
\]
and,
\[
{\hat{\mathcal L}}(s,u,\mathbf X):\mathcal L(s,u,\mathbf X)\left\{\left[1-\bar\zeta(s,Z)\right]\tilde Z(s)+\bar\zeta(s,Z) Z(s)\right\}.
\]
Therefore, $\bar\zeta$ is continuous which implies ${\hat{\mathcal L}}(s,u,\mathbf X)$ is continuous. Since the codomain of ${\hat{\mathcal L}}(s,u,\mathbf X)$ is in $\widetilde\Xi$, every fixed-point of ${\hat{\mathcal L}}(s,u,\mathbf X)$ is also a fixed-point of $\mathcal L$. Now by \emph{$\mathbf{Case\  I}$}, ${\hat{\mathcal L}}(s,u,\mathbf X)$ has a fixed-point and therefore, $\mathcal L$ also does. $\square$

\medskip

\subsection*{Proof of Theorem \ref{t0}}

\medskip

From quantum Lagrangian function expressed in the Equation (\ref{9}), the Euclidean action function for the agent in continuous time $[0,t]$ is given by
\begin{equation*}
\mathcal A_{0,t}(\mathbf X)=\int_0^t\E_s\left\{c\left[u(s),\mathbf X(s)\right]ds+\lambda\left[\bm \mu(s,u,\mathbf X)ds+\bm\sigma(s,\mathbf X) d\mathbf B(s)-\Delta\mathbf Xds\right]\right\},
\end{equation*}
where vector $\lambda>0$ is a time independent quantum Lagrangian multiplier. As at the beginning of the continuous time interval $[s,s+\epsilon]$, as the agent does not have any prior future knowledge, they make expectations based on their  all current state variables represented by $\mathbf X$. Hence, $\E_s[.]:=\E[.|\mathbf X(s),\mathcal F_s]$, where $\mathcal F_s$ is the filtration process starting at time $s$. For a penalization constant $L_\varepsilon>0$ and for time interval $[s,s+\varepsilon]$ with $\varepsilon\downarrow 0$ define a transition function from $s$ to $s+\varepsilon$ as
\begin{equation}\label{13}
\Psi_{s,s+\varepsilon}(\mathbf X)=\frac{1}{L_\varepsilon}\int_{\mathbb R^4}\exp[-\varepsilon\mathcal{A}_{s,s+\varepsilon}(\mathbf X)]\Psi_s(\mathbf X)d\mathbf X,
\end{equation}
where $\Psi_s(\mathbf X)$ is the value of the transition function at time $s$ with the initial condition $\Psi_0(\mathbf X)=\Psi_0$ and the action function in $[s,s+\varepsilon]$ of the representative agent is, 
\begin{equation*}
\mathcal{A}_{s,s+\varepsilon}(\mathbf X)=\int_s^{s+\varepsilon}\E_\nu\left\{c[u(\nu),\mathbf X(\nu)]d\nu+g[\nu+\Delta\nu,\mathbf X(\nu)+\Delta\mathbf X(\nu)]\right\},
\end{equation*}
where $g(\mathbf X)\in C^2([0,t]\times\mathbb{R}^{4})$ such that Assumptions \ref{a0}- \ref{a2} hold and $\widetilde Y(\nu)=g(\mathbf X)$, where $\widetilde Y$ is an It\^o process \citep{oksendal2003} and,
\[
g(\mathbf X)=\lambda \left[\bm \mu(s,u,\mathbf X)ds+\bm\sigma(s,\mathbf X) d\mathbf B(s)-\Delta\mathbf Xds\right]+o(1),
\]
where $\Delta \mathbf X=\mathbf X(s+\varepsilon)-\mathbf X(s)$. In Equation (\ref{13}) $L_\varepsilon$ is a positive penalization constant such that the value of $\Psi_{s,s+\varepsilon}^k(.)$ becomes $1$. One can think this transition function $\Psi_{s,s+\varepsilon}(.)$ as some transition probability function on Euclidean space. We have divided the time interval $[0,t]$ into $n$ small equal sub-intervals $[s,s+\varepsilon]$ so that $\tau=s+\varepsilon$. Fubini's Theorem implies,
\[
\mathcal{A}_{s,\tau}(\mathbf X)=\E_s\left\{\int_s^\tau c[u(\nu),\mathbf X(\nu)]d\nu+g[\nu+\Delta\nu,\mathbf X(\nu)+\Delta\mathbf X(\nu)]\right\}.
\]
After using the fact that $[\Delta\mathbf X(s)]^2=\varepsilon$, for $\varepsilon\downarrow 0$ (with initial condition $\mathbf X(0)$), It\^o's formula and \cite{baaquie1997} imply,
\begin{equation*}
\mathcal{A}_{s,\tau}(\mathbf X)=c[u(s),\mathbf X(s)]+g+\frac{\partial g}{\partial s}+\mu[s,u(s),\mathbf X(s)]\frac{\partial g}{\partial\mathbf X}+\mbox{$\frac{1}{2}$}\bm\sigma^T[s,\mathbf X(s)]\frac{\partial^2 g}{\partial \mathbf X^T\partial\mathbf X}\bm\sigma[s,\mathbf X(s)]+o(1),
\end{equation*}
where $g=g[s,\mathbf X(s)]$. Result in Equation(\ref{13}) implies,
\begin{multline*}
\Psi_{s,\tau}(\mathbf X)=\frac{1}{L_\varepsilon}\int_{\mathbb R^4}\exp\left\{-\varepsilon\left[c[u(s),\mathbf X(s)]+g+\frac{\partial g}{\partial s}+\mu[s,u(s),\mathbf X(s)]\frac{\partial g}{\partial\mathbf X}\right.\right.\\\left.\left.+\mbox{$\frac{1}{2}$}\bm\sigma^T[s,\mathbf X(s)]\frac{\partial^2 g}{\partial \mathbf X^T\partial\mathbf X}\bm\sigma[s,\mathbf X(s)]\right]\right\}\Psi_s(\mathbf X)d\mathbf X+o(\varepsilon ^{1/2}).
\end{multline*}
For $\varepsilon \downarrow 0$ define a new transition probability $\Psi_s^{\tau}$ centered around time $\tau$. A Taylor series expansion (up to second order) of the left hand side of the above Equation yields,
\begin{multline*}
\Psi_s^{\tau}(\mathbf X)+\varepsilon\frac{\partial}{\partial s}\Psi_s^{\tau}+o(\varepsilon)=\frac{1}{L_\varepsilon}\int_{\mathbb R^4}\exp\left\{-\varepsilon\left[c[u(s),\mathbf X(s)]+g+\frac{\partial g}{\partial s}+\mu[s,u(s),\mathbf X(s)]\frac{\partial g}{\partial\mathbf X}\right.\right.\\\left.\left.+\mbox{$\frac{1}{2}$}\bm\sigma^T[s,\mathbf X(s)]\frac{\partial^2 g}{\partial \mathbf X^T\partial\mathbf X}\bm\sigma[s,\mathbf X(s)]\right]\right\}\Psi_s(\mathbf X)d\mathbf X+o(\varepsilon ^{1/2}),
\end{multline*}
as $\varepsilon\downarrow 0$. For fixed $s$ and $\tau$ let $\mathbf X(s)=\mathbf X(\tau)+\vartheta$. For some number $\bar\vartheta^*\in(0,\infty)$ assume $|\vartheta|\leq\bar\vartheta^*\varepsilon[\mathbf X(s)]^{-1}$. Therefore, we get upper bound of state variables in this SIR model as $\mathbf X(s)\leq\bar\vartheta^*\varepsilon/(\vartheta)^2$. Moreover, Fr\"ohlich's Reconstruction Theorem \citep{pramanik2020optimization,pramanik2021optparam,simon1979} and Assumptions \ref{a0}-\ref{a2} imply,
\begin{multline}\label{14}
\Psi_s^{\tau}(\mathbf X)+\varepsilon\frac{\partial}{\partial s}\Psi_s^{\tau}+o(\varepsilon)=\frac{1}{L_\varepsilon}\int_{\mathbb R^4}\exp\left\{-\varepsilon\left[c[u(s),\mathbf X(s)]+g+\frac{\partial g}{\partial s}+\mu[s,u(s),\mathbf X(s)]\frac{\partial g}{\partial\mathbf X}\right.\right.\\\left.\left.+\mbox{$\frac{1}{2}$}\bm\sigma^T[s,\mathbf X(s)]\frac{\partial^2 g}{\partial \mathbf X^T\partial\mathbf X}\bm\sigma[s,\mathbf X(s)]\right]\right\}\left[\Psi_s^\tau(\mathbf X)+\vartheta\frac{\partial}{\partial \mathbf X}\Psi_s^\tau(\mathbf X)+o(\varepsilon)\right]d\mathbf X+o(\varepsilon ^{1/2}),
\end{multline}
as $\varepsilon\downarrow 0$. Define a $C^2$ function,
\begin{equation*}
\tilde f(s,\bar Z)\triangleq c[u(s),\mathbf X(s)]+g+\frac{\partial g}{\partial s}+\mu[s,u(s),\mathbf X(s)]\frac{\partial g}{\partial\mathbf X}+\mbox{$\frac{1}{2}$}\bm\sigma^T[s,\mathbf X(s)]\frac{\partial^2 g}{\partial \mathbf X^T\partial\mathbf X}\bm\sigma[s,\mathbf X(s)].
\end{equation*}
Plugging in $\tilde f(s,\bar Z)$ into Equation (\ref{14}) yields,
\begin{multline}\label{15}
\Psi_s^{\tau}(\mathbf X)+\varepsilon\frac{\partial}{\partial s}\Psi_s^{\tau}(\mathbf X)+o(\varepsilon)=\frac{1}{L_\varepsilon}\int_{\mathbb R^4}\exp\left\{-\varepsilon\tilde f(s,\bar Z)\right\}\\
\times\left[\Psi_s^\tau(\mathbf X)+\vartheta\frac{\partial}{\partial \mathbf X}\Psi_s^\tau(\mathbf X)+o(\varepsilon)\right]d\mathbf X+o(\varepsilon ^{1/2}).
\end{multline}
Let $\tilde f(s,\bar Z)$ be a $C^2$ function. A second order Taylor series expansion yields,
\begin{multline*}
\tilde f(s,u,\vartheta(\tau))=\tilde f(s,u,\vartheta(\tau))+[\vartheta-\mathbf X(\tau)]\frac{\partial}{\partial\mathbf X}\tilde f(s,u,\vartheta(\tau))\\
+\mbox{$\frac{1}{2}$}[\vartheta-\mathbf X(\tau)]^T\frac{\partial^2}{\partial\mathbf X^T\partial\mathbf X}\tilde f(s,u,\vartheta(\tau))[\vartheta-\mathbf X(\tau)]+o(\varepsilon),
\end{multline*}
as $\varepsilon\downarrow 0$ and $\Delta u(s)\downarrow 0$. Define $\hat\vartheta=\vartheta-\mathbf X$ so that $d\hat\vartheta=d\vartheta$. Thus, first integration of Equation (\ref{15}) becomes,
\begin{multline*}
\int_{\mathbf R^4}\exp(-\varepsilon\tilde f(s,\bar Z))d\mathbf X=\int_{\mathbf R^4}\exp\left\{-\varepsilon\left[\tilde f(s,u,\vartheta(\tau))-\mathbf J^T\hat\vartheta+\hat\vartheta^T\mathcal H_{\mathbf X}\hat\vartheta\right]\right\}d\hat\vartheta\\
=\exp\{-\varepsilon\tilde f(s,u,\vartheta(\tau))\}\int_{\mathbb R^4}\exp\left\{(\varepsilon \mathbf J^T)\hat\vartheta-\hat\vartheta(\varepsilon\mathcal H_{\mathbf X})\hat\vartheta\right\}d\hat\vartheta\\=\frac{\pi}{\sqrt{\varepsilon |\mathcal H_{\mathbf X}|}}\exp\left\{\frac{\varepsilon}{4}\mathbf J^T\mathcal H_{\mathbf X}^{-1}\mathbf J-\varepsilon \tilde f(s,u,\vartheta(\tau))\right\},
\end{multline*}
where $\mathbf J=-\partial\tilde f/\partial\mathbf X$ and $\mathcal H_{\mathbf X}$ is a non-singular Hessian matrix. Therefore, first integral term of Equation (\ref{15}) becomes,
\begin{equation*}
\frac{1}{L_\varepsilon}\Psi_s^\tau(\mathbf X)\int_{\mathbb R^4}\exp(-\varepsilon\tilde f)d\mathbf X=\frac{1}{L_\varepsilon}\Psi_s^\tau \frac{\pi}{\sqrt{\varepsilon |\mathcal H_{\mathbf X}|}}\exp\left\{\frac{\varepsilon}{4}\mathbf J^T\mathcal H_{\mathbf X}^{-1}\mathbf J-\varepsilon \tilde f(s,u,\vartheta(\tau))\right\},
\end{equation*}
where $\mathcal H_{\mathbf X}>0$. In a similar fashion we get the second integral term of Equation (\ref{15}) as
\begin{equation*}
\frac{1}{L_\varepsilon}\frac{\partial\Psi_s^\tau(\mathbf X)}{\partial \mathbf X}\int_{\mathbb R^4}\vartheta\exp(-\varepsilon\tilde f)d\mathbf X =\frac{1}{L_\varepsilon}\frac{\partial\Psi_s^\tau}{\partial \mathbf X}\frac{\pi}{\sqrt{\varepsilon |\mathcal H_{\mathbf X}|}}\left[\frac{1}{2}\mathcal H_{\mathbf X}^{-1}+\mathbf X\right]\exp\left\{\frac{\varepsilon}{4}\mathbf J^T\mathcal H_{\mathbf X}^{-1}\mathbf J-\varepsilon \tilde f(s,u,\vartheta(\tau))\right\}.
\end{equation*}
Using above results and Equation (\ref{15}) we obtain a Fokker-Plank type equation as,
\begin{multline*}
\Psi_s^{\tau}(\mathbf X)+\varepsilon\frac{\partial}{\partial s}\Psi_s^{\tau}(\mathbf X)+o(\varepsilon)=\frac{1}{L_\varepsilon}\frac{\pi}{\sqrt{\varepsilon |\mathcal H_{\mathbf X}|}}\exp\left\{\frac{\varepsilon}{4}\mathbf J^T\mathcal H_{\mathbf X}^{-1}\mathbf J-\varepsilon \tilde f(s,u,\vartheta(\tau))\right\}\\
\times\left\{\Psi_s^\tau(\mathbf X)+\left[\frac{1}{2}\mathcal H_{\mathbf X}^{-1}+\mathbf X\right]\frac{\partial\Psi_s^\tau}{\partial \mathbf X}\right\}+o(\varepsilon ^{1/2}),
\end{multline*}
as $\varepsilon\downarrow 0$. Assuming $L_\varepsilon=\pi/\sqrt{\varepsilon | \mathcal H_{\mathbf X}|}>0$ yields,
\begin{multline*}
\Psi_s^{\tau}(\mathbf X)+\varepsilon\frac{\partial}{\partial s}\Psi_s^{\tau}(\mathbf X)+o(\varepsilon)\\=\left\{1+\frac{\varepsilon}{4}\mathbf J^T\mathcal H_{\mathbf X}^{-1}\mathbf J-\varepsilon \tilde f(s,u,\vartheta(\tau))\right\}\left\{\Psi_s^\tau(\mathbf X)+\left[\frac{1}{2}\mathcal H_{\mathbf X}^{-1}+\mathbf X\right]\frac{\partial\Psi_s^\tau}{\partial \mathbf X}\right\}+o(\varepsilon ^{1/2}),
\end{multline*}
as $\varepsilon\downarrow 0$. Since $\mathbf X\leq\bar\vartheta^*\varepsilon/(\vartheta)^2$ assume $|\mathcal H_{\mathbf X}^{-1}|\leq 2\bar\vartheta^*\varepsilon(1-\vartheta^{-1})$ such that $|(2\mathcal H_{\mathbf X})^{-1}+\mathbf X|\leq \bar\vartheta^*\varepsilon$. Therefore, $|\mathcal H_{\mathbf X}^{-1}|\leq 2\varepsilon \bar\vartheta^* (1-\vartheta^{-1})$ so that $|(2\mathcal H_{\mathbf X}^{-1})+\mathbf X|\downarrow 0$. Hence,
\begin{equation*}
\Psi_s^{\tau}(\mathbf X)+\varepsilon\frac{\partial}{\partial s}\Psi_s^{\tau}(\mathbf X)+o(\varepsilon)=(1-\varepsilon)+\Psi_s^\tau +o(\varepsilon ^{1/2}).
\end{equation*}
The Fokker-Plank type equation of stochastic SIR model with infection dynamics is,
\begin{equation*}
\frac{\partial}{\partial s}\Psi_s^{\tau}(\mathbf X)=-\tilde f(s,u,\vartheta(\tau)) \Psi_s^{\tau}(\mathbf X).
\end{equation*}
The solution of 
\begin{equation}\label{16}
-\frac{\partial }{\partial u}\tilde f(s,u,\vartheta(\tau))\Psi_s^\tau(\mathbf X)=0,
\end{equation}
is an optimal ``lock-down" intensity and vaccination rate. Since, $\vartheta=\mathbf X(s)-\mathbf X(\tau)$ for all $\varepsilon\downarrow 0$, in Equation (\ref{16}) $\vartheta$ can be replaced by $\mathbf X$. As the transition function $\Psi_s^{\tau}(\mathbf X)$ is a solution of the Equation (\ref{16}), the result follows. $\square$

\section*{Data availability}
Office for National Statistics. Coronavirus (covid-19) infection survey, antibody data for UK: 16 March, 2021 data have been used.

\bibliography{bib}
\end{document}